\documentclass[12pt]{amsart}
\usepackage{graphicx}
\usepackage{amsfonts}
\usepackage{amssymb}
\usepackage{latexsym}
\usepackage{amscd}
\usepackage{lscape}

\newtheorem{thm}{Theorem}[section]
\newtheorem{cor}[thm]{Corollary}
\newtheorem{lem}[thm]{Lemma}
\newtheorem{prop}[thm]{Proposition}

\newtheorem*{fsc}{Fatou-Sullivan Classification}
\newtheorem*{fshb}{Fatou-Shishikura Bound}
\newtheorem*{tet}{Conditional Implication}

\newenvironment{pf}{\proof[\proofname]}{\endproof}
\newenvironment{pf*}[1]{\proof[#1]}{\endproof}
\usepackage{euscript}

\usepackage[OT2,OT1]{fontenc}

\newcommand{\cal}[1]{{\mathcal #1}}

\newcommand{\beq}{\begin{equation}}
\newcommand{\eeq}{\end{equation}}
\newcommand{\eref}[1]{(\ref{#1})}

\newcommand{\ve}{\varepsilon}
\newcommand{\de}{\delta}
\newcommand{\al}{\alpha}
\newcommand{\be}{\beta}
\newcommand{\ga}{\gamma}

\newcommand{\om}{\omega}

\theoremstyle{definition}
\newtheorem{defn}{Definition}[section]

\theoremstyle{remark}
\newtheorem{rem}{Remark}[section]

\renewcommand{\deg}{\operatorname{deg}}
\newcommand{\riem}{\hat{\CC}}

\newcommand{\dist}{\operatorname{dist}}

\renewcommand{\mod}{\operatorname{mod}}
\newcommand{\tl}{\tilde}

\newcommand{\eps}{\epsilon}

\renewcommand{\Re}{\operatorname{Re}}

\numberwithin{equation}{section}
\newcommand{\thmref}[1]{Theorem~\ref{#1}}
\newcommand{\propref}[1]{Proposition~\ref{#1}}
\newcommand{\secref}[1]{\S\ref{#1}}
\newcommand{\lemref}[1]{Lemma~\ref{#1}}
\newcommand{\corref}[1]{Corollary~\ref{#1}}
\newcommand{\figref}[1]{Figure~\ref{#1}}

\newcommand{\cM}{{\cal M}}

\newcommand{\cB}{{\cal B}}

\newcommand{\cP}{{\cal P}}
\newcommand{\cC}{{\cal C}}

\newcommand{\cD}{{\cal D}}

\newcommand{\CC}{{\mathbb C}}
\newcommand{\RR}{{\mathbb R}}
\newcommand{\JJ}{{\mathbb J}}
\newcommand{\KK}{{\mathbb K}}
\newcommand{\TT}{{\mathbb T}}
\newcommand{\ZZ}{{\mathbb Z}}
\newcommand{\NN}{{\mathbb N}}
\newcommand{\DD}{{\mathbb D}}

\newcommand{\QQ}{{\mathbb Q}}

\newcommand{\fJ}{{\mathfrak J}}

\newcommand{\Om}{\Omega}

\newcommand{\ra}{\rightarrow}

\newcommand{\ze}{\zeta}
\begin{document}
\addtolength{\evensidemargin}{-0.7in}
\addtolength{\oddsidemargin}{-0.7in}

\title[Computability of  Julia sets]{Computability of Julia sets}
\author{Mark Braverman, Michael Yampolsky}
\thanks{The first author's research is supported by an NSERC CGS scholarship}
\thanks{The second author's research is supported by NSERC operating grant}
\date{\today}
\begin{abstract}
In this paper we settle most of the open questions on algorithmic computability 
of Julia sets. In particular, we present an algorithm for constructing quadratics
whose Julia sets are uncomputable. We also show that a filled Julia set of a polynomial
is always computable.
\end{abstract}
\maketitle

\tableofcontents

\section{Foreword}
Computable planar compacts can be visualized on a computer screen with
an arbitrarily high magnification. Among all computer-generated pictures
with mathematical content Julia sets of rational mappings occupy, perhaps,
the most prominent position. And yet, as we have shown in \cite{BY}, some
of those sets are uncomputable, and so cannot be visualized. In this paper
we present an account of the results on computability of Julia sets.
We survey what was previously known, mostly from our own work, some of it
also joint with I.~Binder; and present new results which settle most of the
previously open questions.

A reader, unfamiliar with the questions of algorithmic computability and
complexity, particularly as applied to questions in analysis, will find an
introduction in \secref{sec:intro-comp}. An introduction to the relevant
concepts of dynamics of rational mappings awaits in \secref{sec:intro-dyn}.
Having thus set the stage, we derive some preliminary results on computability
of Julia sets in \secref{sec:prelim-comp}. In \secref{sec:filled} we 
present our first new result:

\medskip
\noindent
{\it All filled Julia sets of polynomial mappings are algorithmically computable.}

\medskip
\noindent
This answers in the positive a question put to us by John Milnor. In the following
 \secref{sec:siegel} we discuss negative results. In particular, we show:

\medskip
\noindent
{\it There exist computable complex parameters $c$, such that the Julia set
of $J(z^2+c)$ is not algorithmically computable.}

\medskip
\noindent
Previously, we could only show the existence of parameters for which the Julia set
is uncomputable. Now we can give an algorithm for their construction. The question
whether such a thing is possible would be invariably asked by our colleagues after
our talks on the subject. The general feeling was, perhaps, that the uncomputability
of Julia sets is in some ways connected to lack of computability of $c$'s. As we
see now, this is not the case.

In the final \secref{sec:interpret} we interpret our results, and attempt to describe
a toy model for uncomputable Julia sets. We also try to answer the ``na\"\i ve'' question:
``What would a computer really draw, when $J(z^2+c)$ is uncomputable?''
Much of the discussion in this section is motivated by a problem posed to us
by Michael Shub. 

We now have a good understanding of computable properties of rational Julia sets.
In this paper, however, we completely avoid discussing the related computational complexity
questions. In this area we know comparatively little. In \cite{BBY2}, jointly with I.~Binder,
we have shown that computational complexity of quadratic Julia sets can be arbitrarily high.
Evidently, there is an interplay between dynamical properties of Julia sets, and the hardness of
drawing their picture on a computer screen. However, as the work \cite{Brv2} of the first
author suggests, this connection may not be straightforward: it shows that a class of maps with
``bad'' dynamical properties has poly-time computable Julia sets. Undoubtedly, many similar
surprises still await us in this line of investigation.

\subsection*{Acknowledgements} It is our pleasure to thank our friend and colleague Ilia Binder
for the many useful discussions on computability of Julia sets.
We thank John Milnor for posing the question on computability of filled Julia sets to us.
We are grateful to Michael Shub for formulating a question which has motivated much of the
discussion on the ``shape'' of uncomputable Julia sets in this paper.
We also wish to thank many colleagues who have invariably asked us if parameters for 
uncomputable Julia sets can be produced algorithmically.

\section{Introduction to computability}
\label{sec:intro-comp}

One of the main goals of computability theory is to classify problems 
according to whether or not they can be solved algorithmically. In fact, 
such questions existed before computers. A famous example
is  {\em Hilbert's Tenth Problem}:

\smallskip
{\em ``Given a diophantine equation with any number of unknown quantities and with rational integral numerical coefficients: to devise a process according to which it can be determined by a finite number of operations whether the equation is solvable in rational integers." [Bulletin of the American Mathematical Society 8 (1902), 437-479.]} 
\smallskip

In other words: 

\smallskip
{\em Is it algorithmically possible to determine if a given diophantine equation is solvable?}

\smallskip

It is fairly clear what an affirmative answer would mean in this case -- a method 
to check if an equation has a solution.
Giving a negative answer (which turns out to be the correct one)
requires a more formal definition of ``methods" that can be used in the solution -- 
as one would need to prove that none of these methods work. A universaly accepted model was
introduced in 1936 in a seminal work by Turing \cite{Tur} in the form of a {\em Turing Machine}.

\subsection{Discrete computability and the Turing Machine}

The definition of a Turing Machine (TM) is somewhat technical and can be found 
in all texts on computability eg. \cite{Papad,Sip}. The computational 
power of a Turing Machine is equivalent to that of a RAM computer, and one 
can think of it as a program on such a computer. The program can use 
a finite amount of memory at each stage of the computation, but it can 
always request more, and there is no {\em a-priori} limit on the amount of 
memory that the machine uses. In fact, there is a general belief, usually 
referred to as the Church-Turing thesis, which states that any computation 
performed on a physical device can be simulated using a Turing Machine. 
Just as an ordinary computer program, any Turing Machine admits a finite 
description. 

The definition of a Turing Machines gives a natural way of classifying the 
computability of functions in the {\em discrete} setting, such as functions
acting on the set of naturals $\NN$ or the set of finite binary strings $\{0,1\}^*$. 
Namely, a function $f(x)$ is computable, if there exists a TM which takes $x$ as
an input and outputs the value $f(x)$.

Computable functions are sometimes called {\it recursive}. They include simple functions
such as integer arithmetic operations and lexicographical sorting of strings. 
They also include problems that appear to be difficult in practice, but can 
be solved nonetheless if we are willing to wait sufficiently long. These include, 
for example, finding the prime factorization of an integer and finding the optimal 
strategy in the game of Go. 

On the other hand, there are many functions that are not computable. One 
argument to see this is a simple counting argument: any TM has a finite 
description, and hence there are countably many TMs. On the other hand, 
there are uncountably many functions from $\NN$ to $\NN$, or even from 
$\NN$ to $\{0,1\}$ -- and thus ``most" functions are not computable.
It is much more iteresting to have specific examples of non-computability. 

One such example is the {\em Halting Problem}. The halting function $H$ maps 
a pair $(T,w)$ where $T$ is an encoding of a TM $M$ and $w$ is a binary input to 
$1$ if the machine $M$ running on input $w$ eventually halts, and $0$ otherwise.

\begin{proof}[Sketch of proof that $H$ is not computable] The proof is by a simple diagonalization 
argument. Suppose there were a TM $M_1$ computing the halting 
function. Let $M_2$ be the following machine: on an input $w$, $M_2$ 
uses $M_1$ to compute $H(w,w)$. If $H(w,w)=0$, then $M_2$ halts, otherwise
it goes into an infinite loop. 

Let $w_2$ be the encoding of $M_2$. What will be the outcome of running $M_2(w_2)$?
If $M_2$ halts on $w_2$, then $H(w_2,w_2)=1$, and thus $M_2$ cannot halt on $w_2$
by definition. If $M_2$ fails to halt on $w_2$, then $H(w_2,w_2)=0$, and by its
definition $M_2$ halts on input $w_2$. In either case we arrive at a contradiction.
\end{proof}

Consider a predicate $A: \NN \times \NN \ra \{0,1\}$ defined as follows. On an input 
$(x,t)$, $A$ viewes $x$ as an encoding of a  pair $(M,w)$ of a TM and an input. 
$A(x,t)$ is $1$ if and only if $x$ gives a valid encoding and $M$ halts on $w$ in 
exactly $t$ steps. It is easy to see that $A$ is a computable predicate using a 
simple simulation. On the other hand, computing the predicate
$$
B(x) = \exists t~A(x,t)
$$
is as difficult as solving the Halting Problem, and thus $B$ is non-computable.
This example will be useful later on. More generally, a predicate of the form 
$P(x) = \exists y ~R(x,y)$ for a computable predicate $R(x,y)$ is said to be 
{\em recursively enumerable}. Moreover, $R(x,y)$ can be modified, so that for every $x$
there exists at most one $y$ such that $R(x,y)$ holds.
 We emphasize this by writing $P(x) = \exists! y ~R(x,y)$ 
 Note that any recursive predicate is also recursively 
enumerable.

Another explicit example of a non-conputable function is given by the negative solution 
to Hilbert's Tenth Problem, which is  due to Matiyasevich (see \cite{Mat} for details 
and the history of the problem). 

\begin{thm}
The function that maps an encoding of a diophantine equation $E$ to $1$ if $E$
is solvable and to $0$ otherwise, is non-compuable. 
\end{thm}

One of Turing's original motivations for introducing the Turing Machine was 
classifying real numbers into computable and non-computable ones. A number is 
said to be computable if there exists a TM that writes its (infinite) decimal 
expansion digit by digit. An equivalent, but slightly less representation dependent is the following 
definition. 

\begin{defn}
A real number $\al$ is said to be computable, if there is a computable function
$\phi: \NN \ra \NN$ such that for all $n$, $\left| \al - \frac{\phi(n)}{2^n} \right| < 2^{-n}$.
The set of the computable reals is denoted by $\RR_\cC$. 
\end{defn}
In other words, there exists an algorithm to approximate $\al$ with any desired degree of precision. 
As with discrete functions, ``most" numbers are non-computable, while most ``nice" numbers such as 
$\pi$ and $e$ are. It can be shown that  $\RR_\cC$ with the usual arithmetic operations 
forms a closed real field.

Here we give an extension of the computable numbers that will be useful later on in the 
paper. 

\begin{defn}
A real number $\al$ is said to be {\em right} computable, if there is a computable function 
$\phi: \NN \ra \QQ$ such that 
\begin{itemize}
\item
the sequence $\{\phi(n)\}$ is nonincreasing: $\phi(1)\ge \phi(2) \ge \ldots$; and
\item 
the sequence $\{\phi(n)\}$ converges to $\al$: $\lim_{n\ra \infty} \phi(n) = \al$. 
\end{itemize}
\end{defn}

It is obvious that a computable real number is also right-computable. The 
converse is not true in general:

\begin{prop}
\label{right-comp}
Right computable numbers form a dense subset in  $\RR\setminus\RR_\cC$. 
\end{prop}
\begin{proof}
It is obviously sufficient to present a single right computable number which is
not computable, as then a dense set can be 
produced using simple arithmetic manipulations.
 Let $P(x) = \exists! y ~ R(x,y)$ be a
non-computable predicate on $\NN$ such that $R(x,y)$ is computable, as 
discussed above. Consider the number 
$$
\al = 1 - \sum_{x=1}^\infty P(x) \cdot 4^{-x}. 
$$
Then $\al$ is non-computable, since computing $\al$ would also enable us to compute 
the predicate $P$. On the other hand, $\al$ is right computable, as demonstrated by the 
following computable function:
$$
\phi(n) = 1 - \sum_{x=1}^{n} \sum_{y=1}^{n} R(x,y) \cdot 4^{-x}. 
$$
$\phi(n)$ is obviously non-increasing, and 
$$
\lim_{n \ra \infty} \phi(n) = 1 - \sum_{x=1}^{\infty} \sum_{y=1}^{\infty} R(x,y) \cdot 4^{-x} = 
1 - \sum_{x=1}^{\infty} P(x) \cdot 4^{-x}  = \al. 
$$
\end{proof}

A more detailed discussion on the different extensions of the concept of a computable number
can be found in \cite{Wei}.

The above definition of computability using Turing Machines  directly applies only to
computability questions for discrete objects. It has to be extended 
if we want to discuss computability of continuous objects such as functions 
over $\RR$ or subsets of $\RR^k$. 

\subsection{Oracle computation, computable real functions}
The history of defining computability for real objects probably begins with 
the work of Banach and Mazur \cite{BM} of 1937, only one year after Turing's paper.
This work has founded the tradition of Computable Analysis (sometimes also called 
Constructive Analysis). Interrupted by war, it was further developed in the book
by Mazur \cite{Maz}. Much research took place in the mid 1950's in the works
of Grzegorczyk \cite{Grz}, Lacombe \cite{Lac}, and others. A parallel school of 
Constructive Analysis was founded by A.~A.~Markov in Russia in the late 1940's.
A modern treatment of the field can be found in \cite{Ko} and \cite{Wei}.

The definition of computability over the reals presented here falls into this framework.

Consider the simplest case in which we would like to compute a function $f: \RR \ra \RR$. 
On an input $x$, we are trying to compute $f(x)$. As in the case with real numbers, the 
machine $M$ computing $f$ should be able to output $f(x)$ with any given precision $2^{-n}$. 
The machine $M$, as well as a practical computer, can only handle a finite amount of information, 
and thus is not capable of reading or storing an entire input $x$. Instead, it is allowed to request 
the input $x$ with an arbitrarily high precision. In other words, it has an external tape and 
a command {\em READ(m)} which requests a $2^{-m}$-approximation $\phi(m)$ of $x$ to be written on 
this tape. It can then be read by the machine from the external tape. It is convenient 
to take all the approximations from the {\em dyadic} set $\DD = \{\frac{k}{2^l}~:~k\in \ZZ,~l\in \NN\}$, as
they possess a natural finite binary encoding. 


To formally define 
computability of real functions let us first introduce the notion of an {\it oracle}:

\begin{defn}
A dyadic-valued function $\phi: \NN \ra \DD$ is called an oracle for a real number $x$
if it satisfies $|\phi(m)-x|<2^{-m}$ for all 
$m$.
\end{defn}

\noindent
An {\it oracle Turing Machine} is a TM which can query the value $\phi(m)$ of some oracle 
$\phi$ for an arbitrary $m\in\NN$. Note that the oracle $\phi$ itself is not a part of the
algorithm, but rather enters as a parameter. We will use a notation $M^\phi$ to 
empasize the dependence of the output of the TM on the values of the oracle.

To get used to the terminology,
imagine a trivial algorithm which given an $n\in\NN$ and a  good enough approximation of $x\in \RR$ outputs a $2^{-n}$-approximation of the number $2x$.
The algorithm executes the command
 $$\text{\em READ }x\text{\em  WITH PRECISION }2^{-(n+1)}.$$
At this point the user (playing the role of an oracle in the dictionary sense)
 enters from the keyboard a dyadic rational $d$ 
for which $|d-x|<2^{-(n+1)}.$ The algorithm proceeds to output $2d$ as the answer.

\begin{defn}
\label{funcomp}
Let $S$ be a subset of $\RR$, and let $f:S \ra \RR$ be a real-valued function on $S$. 
Then $f$ is said to be computable if there is an oracle Turing Machine $M^{\phi}(n)$ 
such that the following holds. If $\phi$ is an oracle for $x\in S$, then for every $n\in\NN$
$M^{\phi}(n)$ returns a dyadic number $q$ such that $|q-f(x)|<2^{-n}$. 
\end{defn}

Note that $M^{\phi}$ is supposed to work with {\em any} valid oracle $\phi$ for $x$. 
The definition generalizes trivially to functions with $k>1$ variables. 

Examples of computable functions include most common functions such as an integer power, $\text{exp}(x)$,
and any trigonometric function. A constant function $f(x)\equiv a$ is computable if and only if
$a$ is a computable number. 

The oracle terminology allows us to separate the problem of computing the parameter $x$ from 
the problem of computing the function $f$ on a {\em given} $x$. For example, the function $x\mapsto x^2$ is computable. Hence even if $a$ is a 
non-computable number, we are still able to compute $a^2$, provided we have an oracle access 
to $a$. This is despite the fact that $a^2$ is a non-computable number. 

A fundamental fact about computable functions in this setting is that computable functions are 
continuous:

\begin{thm}
\label{thm:cont}
Let $S\subset \RR^k$, and suppose $f: S \ra \RR$ is computable by an oracle machine $M^{\phi}$. 
Then $f$ is continuous on $S$. 
\end{thm}

\begin{proof}
Let $x\in S$ and $\ve>0$ be given. Choose an integer $m$ such that $2^{-m}<\ve/2$. 
Let $\phi(n)$ be an oracle for $x$ such that $|\phi(n)-x|<2^{-(n+1)}$ for 
all $n$ (thus ``exceeding" the minimum requirement from an oracle). Then $M^{\phi}(m)$ 
 otputs a number $d \in \DD$ such that $|d-f(x)|<2^{-m}$. It terminates after finitely
 many steps, and hence $\phi$ is only queried up to some finite precision $2^{-k}$. 
 It is now not hard to see that for any $x'$ such that $|x-x'|<2^{-k-1}$, there is 
 a valid oracle $\phi'$ which agrees with $\phi$ up to precision $2^{-k}$. 
 Thus for any $x'\in S \cap (x-2^{-k-1},x+2^{-k-1})$, $M^{\phi'}(m)$ outputs the 
 same answer $d$, and we must have $|d-f(x')|<2^{-m}$. 
 Hence for every $x'\in S$ such that $|x-x'|<2^{-k-1}$, we have
 $$
 |f(x)-f(x')| \le | f(x)-d| + |f(x')-d| < 2^{-m} + 2^{-m} < \ve. 
 $$
\end{proof}

\noindent
In particular, it shows that discontinuous functions, such as $\arcsin$ or $\chi_\QQ$ cannot 
be computed by a single machine on the whole domain of definition.

\noindent
Same considerations can be used to prove a stronger result:

\begin{thm}
In the conditions of \thmref{thm:cont}
there exists a computable function $\mu(x,k):S\times \NN\ra\NN$
 such that
 $$|f(y)-f(x)|<2^{-k}\text{ whenever }
y\in S\text{ and }|y-x|<2^{-\mu(x,k)}.$$
\end{thm}

\noindent
We will refer to the this property by saying that $f$ has a computable local modulus of 
continuity.

\begin{rem}
In some cases, for example when $S=[0,1]$, the {\em global} modulus of continuity
(or simply the modulus of continuity) 
of $f$ on $S$ is also computable. That is, we can compute a function $\mu:\NN\ra\NN$
such that 
\beq
\label{mod cont}
\text{ for any } 
x,y\in S\text{ with }|x-y|<2^{-\mu(k)}\Rightarrow |f(x)-f(y)|<2^{-k}. 
\eeq
More generally, this is true whenever $S$ is a compact computable set (as will be 
defined in the next section). In particular, this is true whenever $S=[a,b]$ with 
computable endpoints $a$ and $b$, or when $S$ is the unit circle in $\RR^2$. 
\end{rem}

\subsection{Computability of subsets of $\RR^k$}

Let $K \subset \RR^k$ be a compact set. We would like to give a definition
for $K$ being computable. In the discrete case the distinction between computability 
of functions and sets is not as important, since a set $S$ is usually said to be 
computable, or decidable, if and only if its characteristic function $\chi_S$ is computable. 
The same definition would not work over $\RR$, since only continuous functions  
 can be computable, hence $\chi_K$ would not be computable unless $K = \emptyset$. 

We say that a TM M computes the set $K$ if it approximates $K$ in the {\it 
Hausdorff metric}. Recall that the Hausdorff metric is a metric on 
compact subsets of $\RR^k$ defined by 
$$d_H ( X, Y) =  \inf \{\epsilon > 0 | X \subset U_{\epsilon} 
(Y)~~\mbox{and}~~  Y \subset U_{\epsilon}(X)\}.$$
 We approximate $K$ using 
a class $\cC$ of sets which is dense in metric $d_H$ among compact 
sets, and such that elements of $\cC$ have a natural binary encoding. 
Namely $\cC$ is the set of finite unions of dyadic balls:
$$
\cC= \left\{ \bigcup_{i=1}^n \overline{B(d_i, r_i)}~|~~\mbox{where}~~d_i\in\DD^k, 
r_i \in \DD \right\}.
$$
Members of $\cC$ can be encoded as binary strings in a natural way. 
The following definition is equivalent to the set computability definition
given in \cite{Wei}, and in earlier works (e.g. \cite{WeiPaper}).

\begin{defn}
\label{setcomp}
We say that a compact set $K \subset \RR^k$ is computable, if exists a TM 
$M(m)$, such that on input $m$, $M(m)$ outputs an encoding of $C_m \in 
\cC$ such that $d_H (K, C_m) < 2^{-m}$. 
\end{defn}

To illustrate the robustness of this definition we present the following 
two equivalent characterizations of computable sets (see e.g. \cite{thesis}). The first one 
relates the definition to computer graphics. It is made more precise 
in the discussion below. The second one relates the 
computability of sets to the computability of functions as per Definition 
\ref{funcomp}.

\begin{thm}
For a compact $K \subset \RR^k$ the following are equivalent:

(1) $K$ is computable as per definition \ref{setcomp},

(2) (in the case $k=2$) $K$ can be drawn on a computer screen 
with arbitrarily high resolution,

(3) the {\em distance function} $d_K (x) = \inf \{ |x-y|~~|~~y\in K \}$ 
is computable as per definition \ref{funcomp}.
\end{thm}

\begin{figure}[ht]
\begin{center}
\includegraphics[angle=0, scale=0.5]{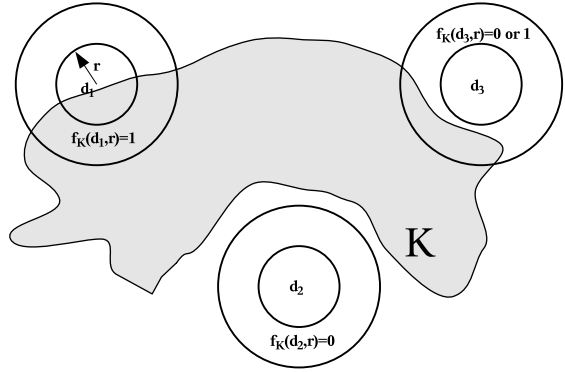}
\end{center}
\caption{Sample values of the function $f_K$}
\label{fig:pixels2}
\end{figure}

Let us elaborate further on part (2) of the theorem. A ``drawing" $P$ of the 
set $K$ on the computer screen is just a collection of pixels that serve 
as an accurate description of $K$ (or a  portion of $K$, if the image is zoomed-in).
We would expect the following properties from $P$:
\begin{itemize}
\item 
$P$ should include all pixels that intersect with $K$, this guarantees
that we get a picture of the entire set $P$; and 
\item 
$P$ should not include pixels that are ``far" from $K$, for example pixels 
that are at least one pixel diameter away from the set $K$.
\end{itemize}
By switching from the rectangular computer pixels to the 
mathematically more convenient round pixels, we see that 
to ``draw" $K$ one should be able to compute a function 
$f_K:\DD \times \DD^k \ra \{0,1\}$ from the family
\begin{equation}
\label{comp:fam}
f_K(d,r) = 
\left\{
\begin{array}{ll}
1 & \mbox{if $B(d,r)\cap K \neq \emptyset$}\\
0 & \mbox{if $B(d,2 \cdot r)\cap K = \emptyset$}\\
0 \mbox{ or } 1 & \mbox{otherwise}
\end{array}
\right.
\end{equation}
$f_K$ then can be used to decide whether to include a round pixel with center 
$d$ and radius $r$ in $P$. Sample values of the function $f_K$ are illustrated
on Figure \ref{fig:pixels2}.

\subsection{Weakly computable sets}

In this section we present a different definition of set-computability, we call 
{\em weak} computability. If was first introduced by Chou and Ko \cite{CK}. 

\begin{defn}
\label{def:weak}
We say that a set $S$ is {\em weakly} computable if there is an oracle
Turing Machine $M^{\phi}(n)$ such that if $\phi=(\phi_1, \phi_2, \ldots, \phi_k)$ represents a point
$x=(x_1, \ldots, x_k) \in \RR^k$, then the output of $M^{\phi}(n)$ is
\beq
\label{weakeq}
M^{\phi}(n) = \left\{
\begin{array}{ll}
1 & \mbox{  if  } x \in K \\
0 & \mbox{  if  } B(x,  2^{-(n-1)}) \cap K = \emptyset \\
0 \mbox{  or  } 1 & \mbox{  otherwise  }
\end{array}
\right.
\eeq
\end{defn}

\begin{figure}[ht]
\begin{center}
\includegraphics[angle=0, scale=0.6]{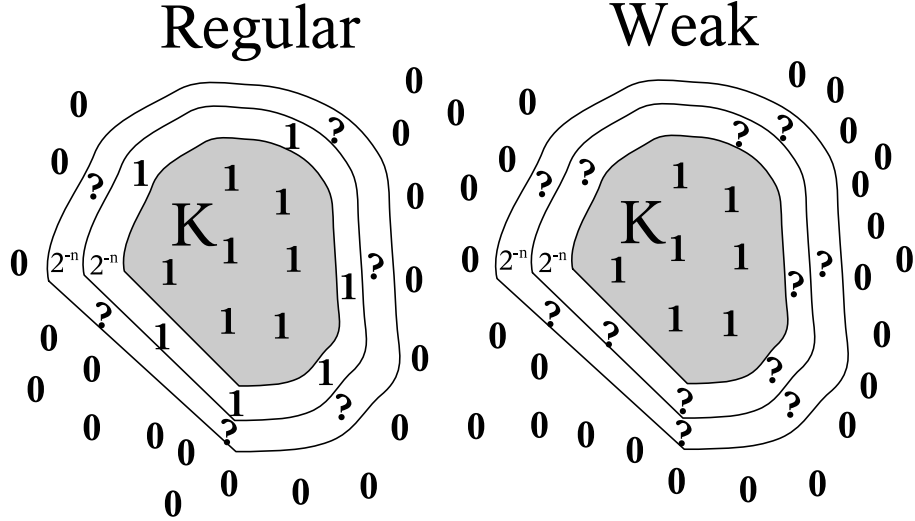}
\end{center}
\caption{The values of $f_K(\bullet,2^{-n})$ in the definitions of regular (left) and weak set computability}
\label{fig:weakstrong}
\end{figure}

Condition \eref{weakeq} is similar to condition \eref{comp:fam}. The 
difference is that now we allow $x$ to be any point in $\RR^k$ (not
just $\DD^k$), and we do not require the machine to output $1$ if 
$x$ is not in $K$ but is ``close". 
It is evident from Figure \ref{fig:weakstrong} that Definition \ref{def:weak} requires
less effort from the algorithm computing $K$ than the original definition.
Thus, the new definition appears to be weaker than the defintion of set computability 
from last section, but it turns out that they are equivalent. 

\begin{thm} {\cite{Brv-FOCS}}
\label{thm:weakstrong}
A compact set $K \subset \RR^k$ is weakly computable if and only if 
it is computable as per definition \ref{setcomp}.
\end{thm}

It is sometimes easier to use weak computability when proving that a 
certain set is computable. We will use Theorem \ref{thm:weakstrong} in 
\secref{sec:interpret}.

\subsection{Set-valued functions and uniformity}

 The problem of computing 
Julia sets is essentially that of mapping the coefficients of a rational 
function $R(z)$ to the set $J_R$. Thus we need a notion of computability 
of set-valued functions to discuss computability questions about Julia sets. 

We can now combine the definitions from previous sections to define 
computability of set valued functions. 

\begin{defn}
\label{def:setval}
Let $S$ be a subset of $\RR^k$. Denote by $K_\ell^*$ the set of all the 
compact subsets of $\RR^\ell$. Let $F: S \ra K_\ell^*$ be a set-valued function
mapping points in $S$ to compact subsets of $\RR^\ell$. $F$ is said to be 
computable on $S$ if there is an oracle TM $M^{\phi_1,\ldots,\phi_k}(n)$ that 
for oracles representing a point $x=(x_1, x_2, \ldots, x_k) \in S$ outputs
an encoding of a set $C_n \in \cC$ such that the Hausdorff distance 
$d_H (F(x),C_n)<2^{-n}$. 
\end{defn}

In fact, the computability definitions for real functions, sets, and set-valued
functions presented above fit in nicely within the much more 
general framework of Type Two Efficiency (TTE). See \cite{Wei}, and references 
therein for more details. In particular, Theorem \ref{thm:cont} stating 
that computable $\Rightarrow$ continuous holds in a very broad variety of 
settings. We will only need it in the case of set-valued functions. The proof
is very similar to the proof of Theorem \ref{thm:cont} (see e.g. \cite{BY}). 

\begin{thm}
\label{cont2}
Suppose $F: S\subset \RR^k \ra K_\ell^*$ is computable as per Definition \ref{def:setval}, 
then $F$ is continous on $S$ in the Hausdorff metric. 
\end{thm}

\noindent
{\bf Example. } Let the complex plane $\CC$ be naturally identified 
with $\RR^2$. Let $d>1$ be an integer. Consider the multi-valued function
$f_d = \sqrt[d]{~}:\CC\ra \CC$. There is no continuous single-valued branch 
of $f_d$ on the entire complex plane, hence there is no computable branch
if $f_d$ that is defined on the entire $\CC$. There is a computable 
branch of $f_d$ that is defined everywhere except for a slit connecting 
$0$ to $\infty$.

On the other hand, if we view the function $f_d$ as a set-valued function
that maps a number $z=r \cdot e^{2 \pi i \theta}$ to its $d$ roots 
$\{ r^{1/d} \cdot e^{2 \pi i \theta/d}, r^{1/d} \cdot e^{2 \pi i (\theta+1)/d}, \ldots,
r^{1/d} \cdot e^{2 \pi i (\theta+d-1)/d}\}$, then it is not hard to 
see that $f_d$ becomes computable. And indeed, the map
$f_d : \RR^2 \ra K_2^*$ is continuous in the Hausdorff metric. 

\smallskip
Note that Definition \ref{def:setval} makes sense even when $S=\{s\}$ is 
a singleton. In this case we say that $F$ is {\em nonuniformly} computable on
$s$. Otherwise, we say that $F$ is {\em uniformly} computable on the set $S$. 

\subsection*{Remark on the BSS computability model} 
We note that another approach to computability of subsets of $\RR^k$
has been developed by Blum, Shub, and Smale \cite{BCSS}. It is based on the
concept of decidability in the Blum-Shub-Smale (BSS) model of real computation.
The BSS model is very different from the Computable Analysis model we use, and
can be very roughly described as based on computation with infinite-precision real arithmetic.
Some discussion of the differences between the models may be found in \cite{BrC} and \cite{Brv-FOCS}.
Algebraic in nature, BSS decidability is not well-suited for the study of 
fractal objects, such as Julia sets. 
It turns out (see Chapter 2.4 of \cite{BCSS}) that in the BSS model all, 
but the most trivial Julia sets 
are not decidable. More generally, sets with a fractional Hausdorff dimension, 
including ones with very simple description, such as the Cantor set, are 
BSS-undecidable.

\section{Julia sets of rational mappings}
\label{sec:intro-dyn}

\subsection{Basic properties of Julia sets}
An excellent general reference for the material in this section is 
the book of Milnor \cite{Mil}.
For a rational mapping $R$ of degree $\deg R=d\geq 2$ considered as a dynamical system on the Riemann
sphere
$$R:\riem\to\riem$$
the Julia set is defined as the complement of the set where the dynamics is Lyapunov-stable:

\begin{defn}
Denote $F(R)$ the set of points $z\in\riem$ having an open neighborhood $U(z)$ on which the
family of iterates $R^n|_{U(z)}$ is equicontinuous. The set $F(R)$ is called the Fatou set of $R$
and its complement $J(R)=\riem\setminus F(R)$ is the Julia set.
\end{defn}

\noindent
In the case when the rational mapping is a polynomial $$P(z)=a_0+a_1z+\cdots+a_dz^d:\CC\to\CC$$ an equivalent
way of defining the Julia set is as follows. Obviously, there exists a neighborhood of $\infty$ on $\riem$
on which the iterates of $P$ uniformly converge to $\infty$. Denoting $A(\infty)$ the maximal such domain of attraction
of $\infty$ we have $A(\infty)\subset F(R)$. We then have 
$$J(P)=\partial A(\infty).$$
The bounded set $\riem \setminus  A(\infty)$ is called {\it the filled Julia set}, and denoted $K(P)$;
it consists of points whose orbits under $P$ remain bounded:
$$K(P)=\{z\in\riem|\;\sup_n|P^n(z)|<\infty\}.$$

\noindent
For future reference, let us summarize in a proposition below the main properties of Julia sets:

\begin{prop}
\label{properties-Julia}
Let $R:\riem\to\riem$ be a rational function. Then the following properties hold:
\begin{itemize}
\item $J(R)$ is a non-empty compact subset of $\riem$ which is completely
invariant: $R^{-1}(J(R))=J(R)$;
\item $J(R)=J(R^n)$ for all $n\in\NN$;
\item $J(R)$ has no isolated points;
\item if $J(R)$ has non-empty interior, then it is the whole of $\riem$;
\item let $U\subset\riem$ be any open set with $U\cap J(R)\neq \emptyset$. Then there exists $n\in\NN$ such that
$R^n(U)\supset J(R)$;
\item periodic orbits of $R$ are dense in $J(R)$.
\end{itemize}
\end{prop}

\noindent
Let us further comment on the last property. For a periodic point $z_0=R^p(z_0)$
of period $p$ its {\it multiplier} is the quantity $\lambda=\lambda(z_0)=DR^p(z_0)$.
We may speak of the multiplier of a periodic cycle, as it is the same for all points
in the cycle by the Chain Rule. In the case when $|\lambda|\neq 1$, the dynamics
in a sufficiently small neighborhood of the cycle is governed by the Mean
Value Theorem: when $|\lambda|<1$, the cycle is {\it attracting} ({\it super-attracting}
if $\lambda=0$), if $|\lambda|>1$ it is {\it repelling}.
Both in the attracting and repelling cases, the dynamics can be locally linearized:
\begin{equation}
\label{linearization-equation}
\psi(R^p(z))=\lambda\cdot\psi(z)
\end{equation}
where $\psi$ is a conformal mapping of a small neighborhood of $z_0$ to a disk around $0$.

\noindent
In the case when $|\lambda|=1$, so that $\lambda=e^{2\pi i\theta}$, $\theta\in\RR$, 
 the simplest to study is the {\it parabolic case} when $\theta=n/m\in\QQ$, so $\lambda$ 
is a root of unity. In this case $R^p$ is not locally linearizable; it is not hard to see that $z_0\in J(R)$.
The description of the dynamics in a small neighborhood of a parabolic orbit 
will be discussed below in some detail.

 In the complementary situation, two non-vacuous possibilities  are considered:
{\it Cremer case}, when $R^p$ is not linearizable, and {\it Siegel case}, when it is.
In the latter case, the linearizing map $\psi$ from (\ref{linearization-equation}) conjugates
the dynamics of $R^p$ on a neighborhood $U(z_0)$ to the irrational rotation by angle $\theta$
(the {\it rotation angle})
on a disk around the origin. The maximal such neighborhood of $z_0$ is called a {\it Siegel disk}.

A different kind of a rotation domain may occur  only for a non-polynomial rational mapping $R$.
A {\it Herman ring} $A$ is a conformal image $$\nu:\{z\in\CC|\;0<r<|z|<1\}\to A,$$
such that $$R^p\circ\nu(z)=\nu(e^{2\pi i\theta}z),$$
for some $p\in\NN$ and $\theta\in\RR\setminus\QQ$.

\noindent
The term {\it basin} in what follows will describe the set of points whose orbits converge to a 
given  periodic orbit under the iteration of $R$. We will denote $\operatorname{Postcrit}(R)$
the {\it post-critical set of $R$}, defined as
the closure of the union of the orbits of critical points of $R$.
Fatou made the following observation:

\begin{prop}
\label{crit orbit}
Let $p_1,\ldots,p_k$ be a periodic orbit of a rational mapping $R$. If it is either attracting, or parabolic,
then its basin contains a critical point of $R$. 
\end{prop}

\noindent
By a perturbative argument, Fatou then concluded that for a rational mapping $R$ with $\deg R=d\geq 2$
at most finitely
many periodic orbits are non-repelling. A sharp bound on their number  depending on $d$ has
been established by Shishikura; it is equal to the number of critical points of $R$ counted with
multiplicity: 

\begin{fshb}
For a rational mapping of degree $d$ the number of the non-repelling periodic
cycles taken together with the number of cycles of Herman rings is at most $2d-2$.
For a polynomial of degree $d$ the number of non-repelling periodic cycles in $\CC$
is at most $d-1$.
\end{fshb}

\noindent
 Therefore, we may refine the last statement of \propref{properties-Julia}:

\begin{itemize}
\item {\it repelling periodic orbits are dense in $J(R)$}.
\end{itemize}

\noindent
Classical results of Fatou also imply the following:

\begin{prop}
\label{cremer postcrit}
Every Cremer point of a rational mapping $R$ as well as every point of the boundary of a Siegel disk
or a Herman ring is contained in $\operatorname{Postcrit}(R)$.
\end{prop}

By definition, the basin of an attracting or a parabolic point, as well as preimages of Siegel disks
and Herman rings belong to the Fatou set. 
 Fatou-Sullivan Classification Theorem formulated below rules out other possibilities:

\begin{fsc}
For every connected component $W\subset F_R$ there exists $m\in\NN$ such that the image $H=R^m(W)$ is
periodic under the dynamics of $R$. Moreover, each periodic Fatou component $H$ is of one of the following
types: 
\begin{itemize}
\item a component of the basin  of  an attracting or a super-attracting periodic orbit;
\item a component of the  basin  of a parabolic periodic orbit;
\item a Siegel disk;
\item a Herman ring.
\end{itemize}
\end{fsc}

To conclude the discussion of the basic properties of Julia sets, let us consider the simplest
examples of non-linear rational endomorphisms of the Riemann sphere, the quadratic polynomials.
Every affine conjugacy class of quadratic polynomials has a unique representative of the
form $f_c(z)=z^2+c$, the family
$$f_c(z)=z^2+c,\;c\in\CC$$
is often referred to as {\it the quadratic family}.
For a quadratic map the structure of the Julia set is governed by the behavior of the orbit of the only
finite critical point $0$. In particular, the following dichotomy holds:

\begin{prop}
\label{quadratic-Julia}
Let $K=K(f_c)$ denote the filled Julia set of $f_c$, and $J=J(f_c)=\partial K$. Then:
\begin{itemize}
\item $0\in K$ implies that $K$ is a connected, compact subset of the plane with connected complement;
\item $0\notin K$ implies that $K=J$ is a planar Cantor set.
\end{itemize}
\end{prop}

\noindent
The {\it Mandelbrot set} $\cM\subset \CC$ is defined as the set of parameter values $c$ for which 
$J(f_c)$ is connected.

\noindent
A rational mapping $R:\hat\CC\to\hat\CC$ is called {\it hyperbolic} if the orbit of every critical point of $R$
is either periodic, or  converges to an
(super-)attracting cycle. The term ``hyperbolic'' has an established meaning in dynamics. Its use in this
context is justified by the following proposition:

\begin{prop}
A rational mapping $R$ of degree $d\geq 2$  is hyperbolic if and only if 
there exists a smooth metric $\mu$ defined on an open neighborhood of $J(R)$ and constants $C>0$, $\lambda>1$
such that
$$||DR^n(z)||_\mu>C\lambda^n\text{ for every }z\in J(R), n\in\NN.$$  
\end{prop}

\noindent
As easily follows from Implicit Function Theorem and considerations of local dynamics of an attracting orbit,
hyperbolicity is an open property in the parameter space of rational mappings of degree $d\geq 2$.

\noindent
Considered as a rational mapping of the Riemann sphere, a quadratic polynomial $f_c(z)$ has two critical points:
the origin, and the super-attracting fixed point at $\infty$. In the case when $c\notin \cM$, the orbit of the
former converges to the latter, and thus $f_c$ is hyperbolic. 
\propref{crit orbit} implies that whenever $f_{c}$ has an attracting orbit in $\CC$, it is a hyperbolic
mapping and $c\in \cM$. The following conjecture is central to the field of dynamics in one complex variable:

\medskip
\noindent
{\bf Conjecture (Density of Hyperbolicity in the Quadratic Family).} {\it Hyperbolic parameters are dense  
in $\cM$.}

\medskip
\noindent
Fatou-Shishikura Bound implies that a quadratic polynomial has at most one non-repelling cycle in the complex
plane. Therefore, 
we will call the polynomial $f_c$ (the parameter $c$, the Julia set $J_c$) {\it Siegel, Cremer,} or 
{\it parabolic} when it has an orbit of the corresponding type.

\subsection{Local dynamics of a parabolic orbit} 
We will describe here briefly the local dynamics of a rational mapping $R$ with a parabolic periodic
point $p$. By replacing $R$ with its iterate, if needed, we may assume that $R(p)=p$, and $R'(p)=1$.
The map $R$ then can be written as
$$R(z)=z+a(z-p)^{n+1}+O((z-p)^{n+2}),\text{ for some }n\in\NN\text{ and }a\neq 0.$$
Note that the integer $n+1$ is the local {\it multiplicity} of $p$ as the solution of $R(z)=z$.

A complex number $\nu\in\TT$ is called an {\it attracting direction} for $p$ if the product $a\nu^n<0$,
and a {\it repelling direction} if the same product is positive.
For each infinite orbit $\{R^{k}(z)\}$ which converges to the parabolic point, there is one
of the $n$ attracting directions $\nu$  for which the unit vectors 
$$(R^{k}(z)-p)/|R^{k}(z)-p|\underset{k\to\infty}{\longrightarrow}\nu.$$
We say in this case that the orbit converges to $p$ in the direction of $\nu$.
For each  attracting direction $\nu$, we say that 
a topological disk $U$ is  an {\it attracting petal} of $R$ at $p$  if the following properties hold:
\begin{itemize}
\item $\overline{U}\ni\{p\}$;
\item  $R^n(\overline{U})\subset U\cup\{p\}$;
\item  an infinite orbit $\{R^{k}(z)\}$ is eventually contained 
in $U$ if and only if it converges to $p$ in the direction of $\nu$.
\end{itemize}

\noindent
 Similarly, $U$ is a {\it repelling
petal} for $R$ if it is an attracting petal for the local branch of $R^{-1}$ which fixes $p$.

\begin{figure}[ht]
\label{fig-flower}
\centerline{\includegraphics[width=0.5\textwidth]{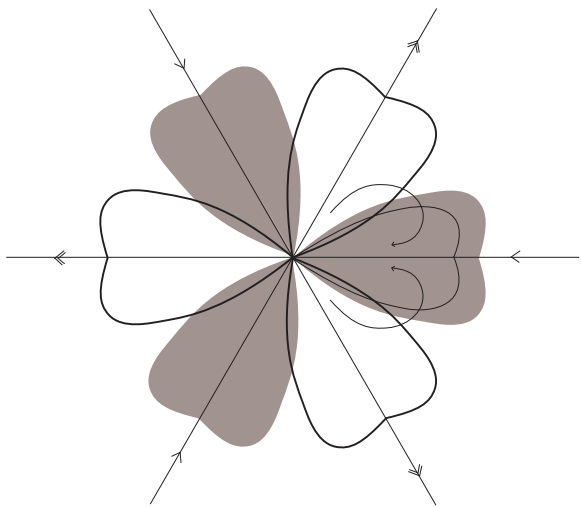}}
\caption{A Leau-Fatou flower with three attracting petals (shaded) and three repelling petals
(emphasized). The attracting and repelling directions are also indicated. The arrows show the direction
of the orbits in one of the petals; the image of this petal is also indicated.}
\end{figure}

The petals form a {\it Leau-Fatou Flower} at $p$:

\begin{thm}
There exists a collection of $n$ attracting petals $P^a_i$, and $n$ repelling petals $P^r_j$
such that the following holds. Any two repelling petals do not intersect, and every repelling petal
intersects exactly two attracting petals. Similar properties hold for attracting petals.
The union
$$(\cup P^a_i)\cup (\cup P^r_j)\cup\{p\}$$
forms an open simply-connected nighborhood of $p$.

\end{thm}

\noindent
The proof of this statement is based on a multivalued change of coordinates
$$w=\kappa(z)=\frac{c}{(z-p)^n},\text{ where }c=-\frac{1}{na}.$$
The map $\kappa$ conformally transforms the infinite sector between two repelling directions into
the plane with the negative real axis removed. In this sector, it changes the map $R$ into
$$F(w)=w+1+O(1/\sqrt[n]{|w|}),\text{ as }w\to\infty.$$
Selecting a right half-plane $H_r=\{\Re z>r\}$ for a sufficiently large $r>0$, we have
$$\Re F(w)>\Re w+1/2,\text{ and hence }F(H)\subset H.$$
The corresponding attracting petal can then be chosen as the domain $\kappa^{-1}(H)$, using the
appropriate branch of the inverse. Note, that given the coefficients of the rational mapping $R$,
the description of the petal is constructive. Let us formulate this last statement in a language
suitable for later references:

\begin{lem}
\label{constructing petal}
For each degree $d\geq 2$ there exists an oracle Turing Machine $M^\phi$ such that the 
following holds. Let $R$ be a rational mapping of degree $d$ with a parabolic periodic point $p$,
with period $m$ and multiplier $e^{2\pi i s/t}$. 
Let $n$ be the number of attracting (and repelling) directions at $p$.
The machine $M^\phi$ takes as input 
the values of $m$, $n$,  $s$, $t$ and a natural number $k$; it is given oracle access to the coefficients
of $R$ and the value of $p$. It outputs a set $L_k\in\cC$ such that the following is true:
\begin{itemize}
\item $L_{k+1}\supset L_k\text{ and }\cup L_k=P\text{ is the union of attracting petals of }R\text{ at }p$,
covering all the attracting directions;
\item $\dist_H(L_k,P)<2^{-k}$.
\end{itemize}
\end{lem}

\noindent
The dynamics inside a petal is described by the following:

\begin{prop}
\label{fatou cyl}
Let $P$ be an attracting or repelling petal of $R$. Then the 
quotient manifold $P/_{z\sim R(z)}$ is conformally isomorphic to the cylinder $\CC/\ZZ$.

\end{prop}

\noindent
In other words, in each of the petals there exists a conformal change of coordinates
transforming $R(z)$ into the unit translation $z\mapsto z+1$.

\medskip
\noindent
Suppose now that the multiplier of the fixed point $p$ is a $q$-th root of unity, $R'(p)=e^{2\pi i p/q}$,
where $(p,q)=1$. A fixed petal for the iterate $R^q$ corresponds to a cycle of $q$ petals for $R$.
It thus follows that $q$ divides the number $n$ of attracting/repelling directions of $p$
 as a fixed point of $R^q$. We make note of
the following proposition, due to Fatou:

\begin{prop}
\label{cycle crit}
Each cycle of attracting petals of a rational mapping $R$ captures an orbit of a critical point of $R$.
\end{prop}

\noindent
This implies, in particular, that a quadratic polynomial $f_c$ with a parabolic periodic point $\zeta$
with multiplier $e^{2\pi i p/q}$ has a Leau-Fatou flower at $\zeta$ with a single cycle of $q$
attracting petals.

\subsection{Occurence of Siegel disks and Cremer points in the quadratic family}
\label{sec:quad siegel}
Let us discuss in more detail the occurrence of Siegel disks in the quadratic family.
For a number $\theta\in [0,1)$ denote $[r_1,r_2,\ldots,r_n,\ldots]$, $r_i\in\NN\cup\{\infty\}$ its possibly finite 
continued fraction expansion:
\begin{equation}
\label{cfrac}
[r_1,r_2,\ldots,r_n,\ldots]\equiv\cfrac{1}{r_1+\cfrac{1}{r_2+\cfrac{1}{\cdots+\cfrac{1}{r_n+\cdots}}}}
\end{equation}
Such an expansion is defined uniquely if and only if $\theta\notin\QQ$. In this case, the {\it rational 
convergents } $p_n/q_n=[r_1,\ldots,r_{n}]$ are the closest rational approximants of $\theta$ among the
numbers with denominators not exceeding $q_n$. In fact, setting $\lambda=e^{2\pi i\theta}$, we have
$$|\lambda^h-1|>|\lambda^{q_n}-1|\text{ for all }0<h<q_{n+1},\; h\neq q_n.$$
The difference $|\lambda^{q_n}-1|$ lies between $2/q_{n+1}$ and $2\pi/q_{n+1}$,
therefore the rate of growth of the denominators $q_n$ describes how well 
$\theta$ may be approximated with rationals.

\begin{defn}
The {\it diophantine numbers
of order k}, denoted $\cD(k)$
is the following class of irrationals ``badly'' approximated by rationals.
 By definition, $\theta\in\cD(k)$ if there exists $c>0$ such that
$$q_{n+1}<cq_n^{k-1}$$
\end{defn}

\noindent
The numbers $q_n$ can
be calculated from the recurrent relation
$$q_{n+1}=r_{n+1}q_n+q_{n-1},\text{ with }q_0=0,\; q_1=1.$$ Therefore, $\theta\in\cD(2)$ if and only if the sequence $\{r_i\}$
is bounded. Dynamicists call such numbers {\it bounded type} (number-theorists prefer {\it constant type}). 
An extreme example of a number of bounded type is the golden mean
$$\theta_*=\frac{\sqrt{5}-1}{2}=[1,1,1,\ldots].$$
The set $$\displaystyle\cD(2+)\equiv\bigcap_{k>2}\cD_k$$
has full measure in the interval $[0,1)$. In 1942 Siegel showed:

\begin{thm}[\cite{siegel}]
Let $R$ be an analytic map with a periodic point $z_0\in\riem$ of period $p$. Suppose the multiplier of the cycle
$$\lambda=e^{2\pi i\theta}\text{ with }\theta\in\cD(2+),$$
then the local linearization equation (\ref{linearization-equation}) holds.
\end{thm}

\noindent
The strongest known generalization of this result was proved by Brjuno in 1972:
\begin{thm}[\cite{Bru}]
Suppose
\begin{equation}
\label{brjuno}
B(\theta)=\displaystyle\sum_n\frac{\log(q_{n+1})}{q_n}<\infty,
\end{equation}
then the conclusion of Siegel's Theorem holds.
\end{thm}

\begin{figure}[h]
\centerline{\includegraphics[width=0.5\textwidth]{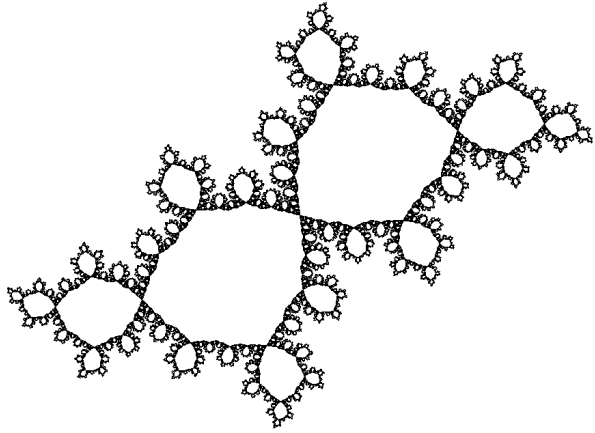}}
\caption{The  Julia set of $P_\theta$ for $\theta=[1,1,1,1,\ldots]$ (the inverse golden mean).}
\label{fig-siegel}
\end{figure}

\noindent
Note that a quadratic polynomial with a fixed Sigel disk with rotation angle $\theta$ after an affine
change of coordinates can be written as 
\begin{equation}
\label{form-1}
P_\theta(z)=z^2+e^{2\pi i \theta}z.
\end{equation}
\noindent
In 1987 Yoccoz \cite{Yoc} proved the following converse to Brjuno's Theorem:

\begin{thm}[\cite{Yoc}]
Suppose that for $\theta\in[0,1)$ the polynomial $P_\theta$ has a Siegel point at the origin.
Then $B(\theta)<\infty$.
\end{thm}

\noindent
The numbers satisfying (\ref{brjuno}) are called Brjuno numbers; the set of all Brjuno numbers will be denoted $\cB$.
It is evident that $\cup\cD(k)\subset \cB$
and thus the set $\cB$ has full measure in the unit circle. On the other hand, it can be shown that its complement
is dense-$G_\delta$.

 The sum of the series (\ref{brjuno}) is called the Brjuno function. 
For us a different characterization of $\cB$ will be more useful. Inductively define $\theta_1=\theta$
and $\theta_{n+1}=\{1/\theta_n\}$. In this way, 
$$\theta_n=[r_{n},r_{n+1},r_{n+2},\ldots].$$
We define the {\it Yoccoz's Brjuno function} as
$$\Phi(\theta)=\displaystyle\sum_{n=1}^{\infty}\theta_1\theta_2\cdots\theta_{n-1}\log\frac{1}{\theta_n}.$$
One can verify that $$B(\theta)<\infty\Leftrightarrow \Phi(\theta)<\infty.$$
The value of the function $\Phi$ is related to the size of the Siegel disk in the following way.

\begin{defn}
Let $P(\theta)$ be a quadratic polynomial with a Siegel disk $\Delta_\theta\ni 0$. Consider a conformal isomorphism
$\phi:\DD\mapsto\Delta$ fixing $0$. The {\it conformal radius of the Siegel disk $\Delta_\theta$} is
the quantity
$$r(\theta)=|\phi'(0)|.$$
For all other $\theta\in[0,\infty)$ we set $r(\theta)=0$. 
\end{defn} 

\noindent
By the Koebe One-Quarter Theorem of classical complex analysis, the internal radius of $\Delta_\theta$ is at least
$r(\theta)/4$. Yoccoz \cite{Yoc} has shown that the sum 
$$\Phi(\theta)+\log r(\theta)$$
is bounded from below independently of $\theta\in\cB$. Recently, Buff and Ch{\'e}ritat have greatly improved this result
by showing that:

\begin{thm}[\cite{BC}]
\label{phi-cont}
The function $\theta\mapsto \Phi(\theta)+\log r(\theta)$ extends to $\RR$ as a 1-periodic continuous
function.
\end{thm}

\noindent
We remark that the following stronger conjecture exists (see \cite{MMY}):

\medskip
\noindent
{\bf Marmi-Moussa-Yoccoz Conjecture.} \cite{MMY} {\it The function $\upsilon: \theta\mapsto \Phi(\theta)+\log r(\theta)$ is H{\"o}lder of exponent $1/2$.}

\medskip
\noindent
Let us remark here, even though we will not use it in the present paper, that in \cite{BY} we have
demonstrated:

\begin{thm}
\label{phi noncomp}
There exists $\theta_0\in\cB$ such that the function $\theta\mapsto\Phi(\theta)$ is uncomputable
on the domain consisting of a single point $\{\theta_0\}$ by a Turing Machine with an oracle 
access to $\theta$.
\end{thm}

\noindent
Assuming Marmi-Moussa-Yoccoz Conjecture holds, \thmref{phi noncomp} would be sufficient to demonstrate that 
$r(\theta)$ is not computable for some values of $\theta\in\TT$; which in turn, by \thmref{noncomp radius} below,
would imply non-computability of $J(P_\theta)$:

\begin{tet}
If the function 
$$\upsilon: \theta\mapsto \Phi(\theta)+\log r(\theta)$$
has a computable modulus of continuity, then it 
is uniformly computable on the entire interval $[0,1]$. 
\end{tet}

\noindent
The proof of the above implication uses the following result of Buff and Ch{\'e}ritat (\cite{BC}). 

\begin{lem}[\cite{BC}]
For any rational point $\theta = \frac{p}{q} \in [0,1]$ denote, as before,
$$P_{\theta}(z)= 
e^{2 \pi i \theta} z + z^2,$$ and let the Taylor expansion of $P_\theta^{\circ q} (z)$ at $0$ start with  
$$
P_\theta^{\circ q} (z) = z + A z^{q+1} + \ldots,\text{ for }q\in\NN
$$
Let $L(\theta) = \left( \frac{1}{q A}\right)^{1/q}$. Denote by $\Phi_{trunc}$
the modification of $\Phi$ applied to rational numbers where the sum is truncated
before the infinite term. Then we have the following explicit formula
for computing $\upsilon(\theta)$:
\beq
\label{expl-rat}
\upsilon(\theta) = \Phi_{trunc} (\theta) + \log L(\theta) + \frac{\log 2\pi}{q}.
\eeq
\end{lem}

\noindent
Equation \eref{expl-rat} allows us to compute the value of $\upsilon$ easily at
every rational $\theta \in \QQ \cap [0,1]$ with an arbitrarily good precision. 
Assuming that $\upsilon$ has a computable modulus of continuity,
it is computable by a single machine of the interval $[0,1]$ 
(see for example Proposition 2.6 in \cite{Kosurv}). This implies the Conditional
Implication.

\noindent
The following conditional result follows:

\begin{lem}
[{\bf Conditional}] \label{2main-equiv}
Suppose the Conditional Implication holds. Let $\theta \in [0,1]$
be such that $\Phi(\theta)$ is finite. Then there is an oracle Turing 
Machine $M^{\phi}_1$ computing $\Phi(\theta)$ with an oracle access to 
$\theta$ if and only if there is an oracle Turing 
Machine $M^{\phi}_2$ computing $r(\theta)$ with an oracle access to 
$\theta$. 
\end{lem}

\begin{proof}
Suppose that $M^{\phi}_1$ computes $\Phi(\theta)$ for some 
$\theta$. Let $M^{\phi}$ be the machine uniformly computing the 
function $\upsilon$. Then we can use $M_1^{\phi}$ and $M^{\phi}$ 
to compute $\log r(\theta) = \upsilon(\theta) - \Phi (\theta)$ with 
an arbitrarily good precision. We can then use this construction 
to give a machine $M_2^{\phi}$ which computes $r(\theta)$. 

The opposite direction is proved analogously. 
\end{proof}

\begin{figure}[ht]
\label{phi figure}
\centerline{
\mbox{\includegraphics[width=0.45\textwidth]{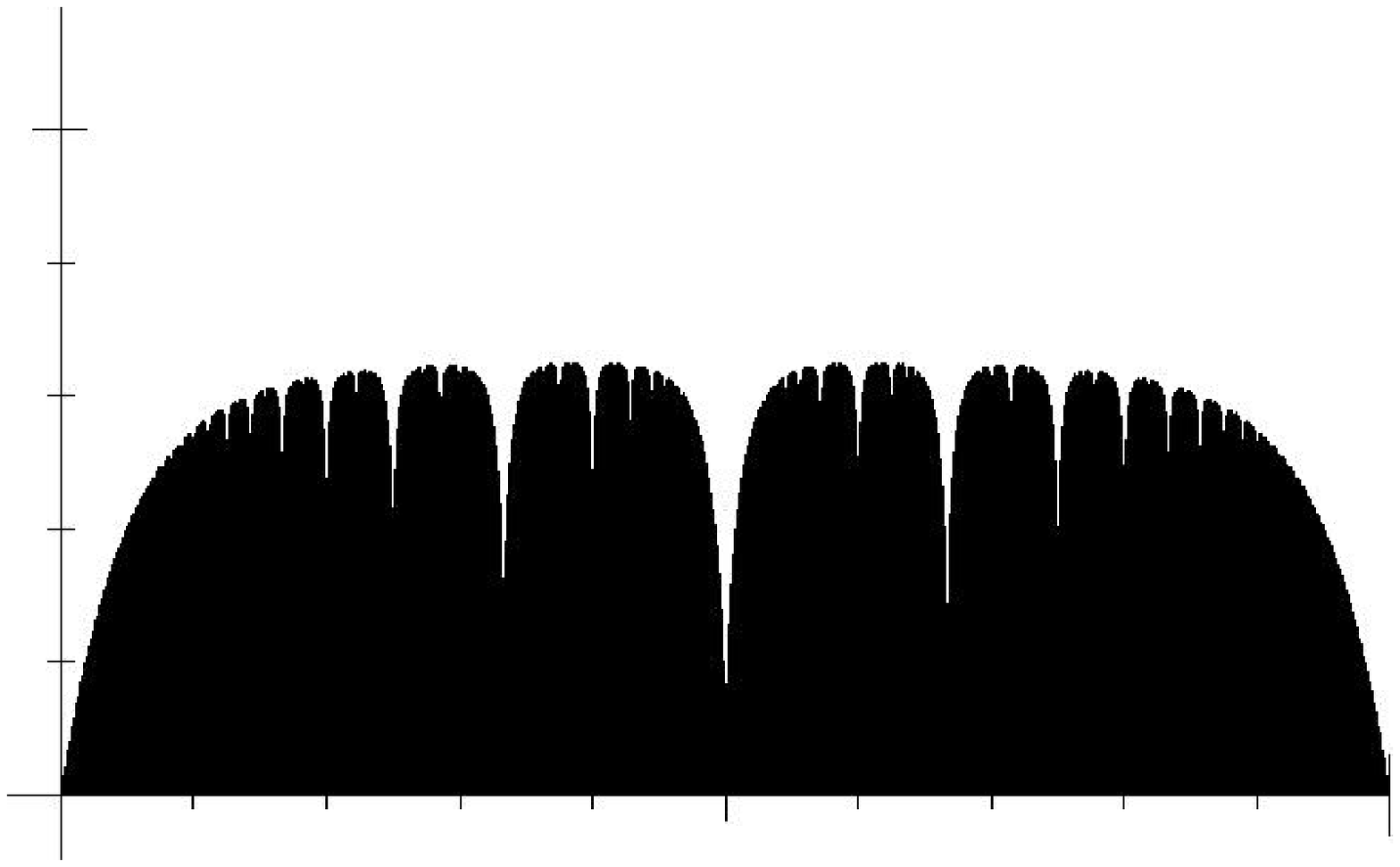}}
\mbox{\includegraphics[width=0.45\textwidth]{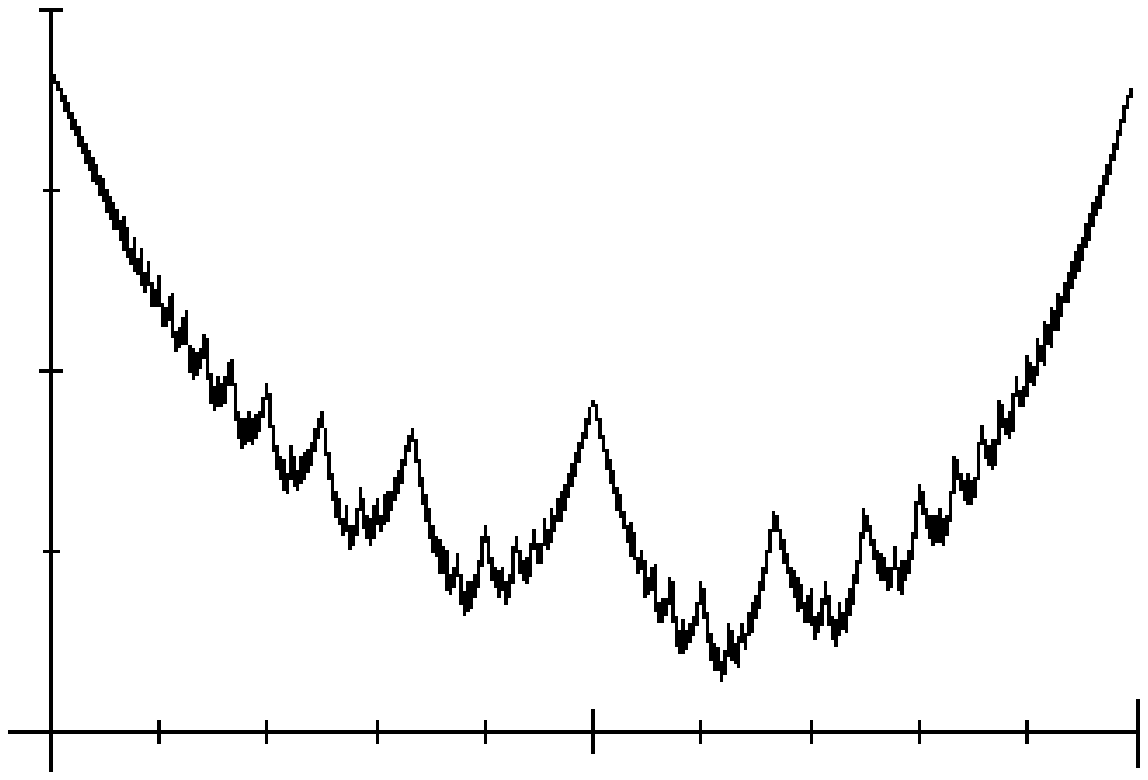}}
}
\caption{The figure on the left is an attempt to visualize the (uncomputable!) function $\Phi$,
by plotting the heights of $\exp(-\Phi(\theta))$ over a grid of Brjuno irrationals.
On the right is the graph of the (conjecturally computable) function $\upsilon(x)$.}

{\it Both figures courtesy of Arnaud Ch{\'e}ritat}

\end{figure}

\subsection*{A note on topological properties of Siegel and Cremer quadratic Julia sets.}
By \propref{cremer postcrit}, the Julia set of any Siegel or Cremer quadratic polynomial is connected.
The following result is due to Sullivan and Douady (see \cite{Sul}):

\begin{thm}
\label{loc con}
If the Julia set of a polynomial mapping $f$ is locally connected, then 
$f$ has no Cremer points. 
Moreover, every cycle of Siegel disks of $f$ contains at least one
critical point in its boundary.

\end{thm} 

\noindent
Thus, in particular, Cremer quadratic Julia sets are {\it never} locally connected. There is a vast
amount of recent work on pathological properties of Cremer quadratics, and we will not attempt to give
a survey of results here. Let us only mention a  paper of S{\o}rensen \cite{Sor} in which
there is a discussion of the mechanism of non local-connectedness in some such sets, which also gives
some indication of the visual complexity of pictures of Cremer Julia sets. We cannot offer an illustration
with a Cremer Julia set to the reader -- even though we will see that all such sets are computable, 
no informative pictures of them have been produced to this day.

As for Siegel Julia sets, Petersen \cite{Pet} showed that $J(P_\theta)$ is locally connected for $\theta$
of bounded type. A different proof of this was later given by the second author \cite{Yam}. Petersen and 
Zakeri \cite{PZ} further extended this result to a set of angles $\theta$ which has a full measure in $\TT$.

On the other hand, Herman in 1986 presented first examples of $P_\theta$ with a Siegel disk whose boundary
does not contain any critical points. By \thmref{loc con} the Julia set of such a map is not locally-connected.
In recent papers of Buff-Ch{\'e}ritat \cite{BC2}, and Avila-Buff-Ch{\'e}ritat \cite{ABC} it is shown that the boundary 
$\partial \Delta_\theta$ of the Siegel disk itself can have smoothness just a hair breadth short of 
analytic in such cases.

For $\theta$ of bounded type, the boundary $\partial\Delta_\theta$ is a quasi-fractal Jordan curve
(a quasi-circle, see \cite{Ahl})
passing through the critical point $c_\theta$ of $P_\theta$. To visualize it, we can use the following
fact:

\begin{prop}
\label{bdary bded type}
Let $\theta$ be of bounded type, and denote $p_n/q_n$ its continued fraction convergents. Let $B>0$
be an upper bound on $\sup q_{n+1}/q_n$. There exist constants $K>0$, $\tau<1$ which depend only 
on $B$, such that
$$\dist_H(\Omega_n,\partial\Delta_\theta)<K\tau^n,\text{ where }\Omega_n=\{P_\theta^i(c_\theta),\;i=0,\ldots,q_{n}\}.$$

\end{prop} 

\noindent
The proof is given in the Appendix.
\propref{bdary bded type} gives a recipe for drawing boundaries of Siegel disks of bounded type,
see, for instance, Figure \ref{fig-siegel}.

\subsection*{Dependence of the conformal radius of a Siegel disk on the parameter}
In this section we will show that the conformal radius of a Siegel disk varies continuously with the Julia set.
To that end we will need a preliminary definition:

\begin{defn}
Let $(U_n,u_n)$ be a sequence of topological disks $U_n\subset\CC$ with marked points $u_n\in U_n$.
The {\it kernel} or {\it Carath{\'e}odory} convergence $(U_n,u_n)\to (U,u)$ means the following:
\begin{itemize}
\item $u_n\to u$;
\item for any compact $K\subset U$ and for all $n$ sufficiently large, $K\subset U_n$;
\item for any open connected set $W\ni u$, if $W\subset U_n$ for infinitely many $n$, then $W\subset U$.
\end{itemize}
\end{defn}

\noindent
The topology on the set of pointed domains which corresponds to the above definition of convergence is again
called {\it kernel} or {\it Carath{\'e}odory} topology. The meaning of this topology is as follows.
For a pointed domain $(U,u)$ denote 
$$\phi_{(U,u)}:\DD\to U$$
the unique conformal isomorphism with $\phi_{(U,u)}(0)=u$, and $(\phi_{(U,u)})'(0)>0$.  
We again denote $r(U,u)=|(\phi_{(U,u)})'(0)|$ the conformal radius of $U$ with respect to $u$.

By the Riemann Mapping
Theorem, the correspondence $$\iota:(U,u)\mapsto \phi_{(U,u)}$$
establishes a bijection between marked topological disks properly contained in $\CC$ and univalent maps $\phi:\DD\to\CC$
with $\phi'(0)>0$.
The following theorem is due to Carath{\'e}odory, a proof may be found in  \cite{Pom}:

\begin{thm}[{\bf Carath{\'e}odory Kernel Theorem}]
The mapping $\iota$ is a homeomorphism with respect to the Carath{\'e}odory topology on domains and the
compact-open topology on maps.
\end{thm}

\noindent
\begin{prop}
\label{radius-continuous}
The conformal radius of a quadratic Siegel disk varies continuously with respect to the Hausdorff 
distance on Julia sets.
\end{prop}

\begin{pf}
To fix the ideas, consider the family $P_\theta$ with $\theta\in\cB$ and denote $\Delta_\theta$ the Siegel disk
of $P_\theta$. It is easy to see that the Hausdorff convergence $J(P_{\theta_n})\to J(P_\theta)$ implies the
Carath{\'e}odory convergence of the pointed domains
$$(\Delta_{\theta_n},0)\to(\Delta,0).$$
The proposition follows from this and the Carath{\'e}odory Kernel Theorem.
\end{pf}

\noindent
In fact, we can state the following quantitative version of the above result. 
For a pointed domain $(U,u)$ denote $\rho(U,u)$ the {\it inner radius}
$\rho(U,u)=\dist(u,\partial U)$.
\begin{lem}
\label{variation conf radius}
Let $U$ be a simply-connected  bounded subdomain of $\CC$ containing the point $0$ in the interior.
Suppose $V\subset U$ is a simply-connected subdomain of $U$, and $\partial V\subset B(\partial U,\eps)$.
Then 
$$r(U,0)-r(V,0)\leq 4\sqrt{r(U,0)}\sqrt{\eps}.$$
Moreover, denote $F(x)=4x/(1+x)^2.$ Then
$$r(V,0)\leq r(U,0)F\left(\frac{\rho(V,0)}{\rho(U,0)} \right).$$

\end{lem}
\noindent
The first inequality is based on Koebe Theorem, see e.g. \cite{RZ} for a proof. The left-hand side is 
a standard refinement of Schwarz Lemma.

\noindent
An immediate corollary is:
\begin{cor}
\label{noncomp radius}
Suppose the function $r(\theta)$ is uncomputable on the set $\{\theta_0\}$. Then the function
$\theta\mapsto J(P_\theta)$ is also uncomputable at the same point.

\end{cor}
\begin{pf}
Assume that $J(P_{\theta_0})$ is computable. Using the output of the TM computing this Julia set
in an obvious way,
for each $\eps>0$ we can obtain a  domain $V\in\cC$ such that 
$$V\subset \Delta_{\theta_0}\text{ and }
d_H(\partial V,\partial \Delta_{\theta_0})<\eps.$$ 
It is elementary to verify that for every $\theta\in\TT$, the set $J(P_\theta)\subset B(0,2)$.
This implies, by Schwarz Lemma, that the conformal radius $r(\theta_0)<2$. Hence, 
by \lemref{variation conf radius}, 
$$|r(V,0)-r(\theta_0)|<\delta=8\sqrt{\eps}.$$
Using any constructive version of the Riemann Mapping Theorem (see e.g. \cite{BB}), we 
can compute $r(V,0)$ to precision $\delta$, and hence know $r(\theta_0)$ up to an error of $2\delta$.
Given that $\delta$ can be made arbitrarily small, we have shown that $r(\theta_0)$ is computable.

\end{pf}

\noindent
We also state for future reference the following proposition:

\begin{prop}
\label{r doest drop}
Let $\{\theta_i\}$ be a sequence of Brjuno numbers such that $\theta_i\to\theta$ and
$\overline{\lim}\; r(\theta_i)=l>0$. Then $\theta$ is also a Brjuno number and $r(\theta)\geq l$.
\end{prop}
\begin{pf}
Denote $\phi_i\equiv \phi_{(\Delta_{\theta_i},0)}.$ 
Note that by Schwarz Lemma, the inverse $\psi_i\equiv (\phi_i)^{-1}$ linearizes $P_{\theta_i}$
on $\Delta_{\theta_i}$. By passing to a subsequence we can assure that 
$\phi_i\to\phi$ locally uniformly, and $\phi'(0)\geq l$. By continuity, $\phi^{-1}$ is a linearizing
coordinate for $P_\theta$, so $\theta$ is a Brjuno number. Moreover, $\phi(\DD)\subset \Delta_\theta$,
and so by Schwarz Lemma $r(\theta)\geq l$.
\end{pf}

\section{Preliminary results on computability of Julia sets}
\label{sec:prelim-comp}

\subsection{Computability without oracle access to $c$}
It is a natural question to ask how easy or how difficult it is to draw a picture of a 
quadratic Julia set {\it without} an oracle access to the value of $c$. As we see below,
in such conditions even very simple Julia sets become algorithmically uncomputable.
Note first the following elementary statement:
\begin{prop}
If $c\in(-\infty,-2)$ then $f_c$ is hyperbolic, and $J_c$ is a Cantor set.
Moreover,
$J_c\subset \overline{B(0,\beta_c)}$, where $\beta_c=\sqrt{1/4-c}+1/2>2$ is a fixed point of $f_c$.
\end{prop}
\begin{pf}
Let $z\in \CC$ with $|z|=\beta_c+\delta$, for some $\delta>0$. By the Triangle Inequality,
$$|f_c(z)|=|z^2+c|\geq |z^2|+|c|=|z|^2+c=(\beta_c+\delta)^2+c>$$
$$>\beta_c^2+c+2\beta_c\delta=f_c(\beta_c)+2\beta_c\delta>\beta_c+4\delta.$$
It follows immediately that $f_c^n(z)\to\infty$, and hence $J_c\subset \overline{B(0,\beta_c)}$.
It remains to note that
 $$c=\beta_c(1-\beta_c)<-\beta_c,\text{ and hence }f_c(c)>\beta_c.$$
\end{pf}

\begin{thm}
\label{uncomput no oracle}
Let $c<-2$ be an uncomputable real number. Then the Julia set $J_c$ is uncomputable 
by a Turing Machine without oracle access to $c$.

\end{thm}

\begin{pf}
The fixed point $\beta_c=\sqrt{1/4-c}+1/2$ of the mapping $f_c$ is repelling
under our assumption on $c$, and hence lies in the Julia set. 
By the previous proposition, 
$$\beta_c=\sup_{z\in J_c}|z|.$$
Now assume that there exists a Turing Machine $M(n)$ which computes $J_c$.
Use it to determine the largest $j>0$ such that $j\cdot 2^{-n}$ is at most
$2^{-n}$-far from all points in $J_c$. Then 
$$0<(j\cdot 2^{-n}-\beta_c)<2^{-(n-1)},$$
hence, $\beta_c$ is computable. But
$$c=\beta_c-\beta_c^2,$$
which contradicts the assumption that $c$ is an uncomputable real.

\end{pf}

\subsection{Lack of uniform computability of Julia sets}
Another natural question to consider is whether it is possible to compute {\it all}
 Julia sets, or in particular all quadratic Julia sets,
 with a single oracle Turing Machine $M^\phi(n)$. This is ruled out by \thmref{cont2}, as
the dependence $c\mapsto J(f_c)$ is discontinuous in the Hausdorff distance.
For an excellent survey of this problem see the paper of Douady \cite{Do}. 

\begin{thm}[\cite{Do}]
\label{discont}
Denote $\JJ(c)$ and $\KK(c)$  the functions $c\mapsto J_c$ and $c\mapsto K_c$ respectively
viewed as functions from $\CC$ to $K^*_2$ with the latter space equipped
with Hausdorff distance. 
Then the following is true:
\begin{itemize}
\item[(a)] if $c$ is Siegel then $\JJ(c)$ is discontinuous at $c$, but $\KK(c)$ is continuous at $c$;
\item[(b)] if $c$ is parabolic then both $\JJ(c)$ and $\KK(c)$ are discontinuous at $c$;
\item[(c)] if $c$ is neither Siegel, nor parabolic, then both $\JJ(c)$ and $\KK(c)$ are continuous at $c$.

\end{itemize}
\end{thm}

\noindent
The discontinuity of $\JJ$ at Siegel parameters is not difficult to prove:

\begin{prop} 
\label{siegel discont}
Let $c_*\in\cM$ be a parameter value for which $f_c$ has a Siegel disk. Then the map $\JJ(c)$
is discontinuous at $c_*$. More specifically, 
let $z_0$ be the center of the Siegel disk.
For each $s>0$ there exists $\tl c\in B(c,s)$ such that
$f_{\tl c}$ has a parabolic periodic point in $B(z_0,s)$.
\end{prop}

\begin{pf}
Denote $\Delta$ the Siegel disk around $\zeta_0$, 
$p$ its period, and $\theta$ the rotation angle.
By the Implicit Function Theorem, there exists a holomorphic mapping $\zeta:U(c_*)\to \CC$ such that
$\zeta(c_*)=z_0$ and $\zeta(c)$ is fixed under $(f_c)^p$. The mapping 
$$\nu:c\mapsto D(f_c)^p(\zeta(c))$$
is holomorphic, hence it is either constant or open. If it is constant, all quadratic polynomials have a Siegel
disk. This is not possible: for instance, $f_{1/4}$ has a parabolic fixed point, and thus no other
non-repelling cycles. Therefore, $\nu$ is open, and in particular, there is a sequence of 
parameters $c_n\to c_*$ such that $\zeta(c_n)$ has multiplier $e^{2\pi i p_n/q_n}$. 
Since $\zeta(c_n)$ is parabolic, it lies in the Julia set of $f_{c_n}$. Hence
$$\dist_H(J(f_{c_*}),J(f_{c_n}))>\dist(c_*,\partial\Delta)/2$$
for $n$ large enough.
\end{pf}

\noindent
Thus an arbitrarily small change of the multiplier of the Siegel point may lead to an implosion of the 
Siegel disk -- its inner radius collapses to zero.

As an immediate consequence of Proposition \ref{siegel discont} and
\thmref{cont2} we have:

\begin{prop}
For any TM $M^\phi(n)$ with an oracle for $c\in\CC$ denote $S_M$ the set
of all values of $c$ for which $M^\phi$ computes $J_c$. Then $S_M\neq \CC$.
\end{prop}

\noindent
In other words, a single algorithm for computing all quadratic Julia sets
does not exist.

\subsection{Discontinuity at a parabolic parameter}
\label{sec:implosion}
The discontinuity in $\JJ(c)$ which occurs at parabolic parameter values has found 
many interesting dynamical implications. The proof is very involved, its outline may
be found in \cite{Do}. It is based on the Douady-Lavaurs theory of {\it parabolic
implosion.} Let us briefly describe its mechanism for the case of a quadratic polynomial
$f_c$. 

Denote $\zeta$ a parabolic periodic point of  
$f_c$ with multiplier $e^{2\pi ip/q}$, and let $m\in\NN$ be its period. 
Let $P_A$ and $P_R$ be an attracting and a repelling petals of $f_c$. Recall that 
by \propref{cycle crit}, the cycle of images $f_c^{jm}(P_A\cup P_R)$, $j=0,\ldots,q-1$
forms a full Leau-Fatou flower at $\zeta$. 

By \propref{fatou cyl}, the quotient
$$C_A=P_A/f_c^{mq}\simeq \CC/\ZZ.$$
The quotient $C_A$, is sometimes called the {\it attracting Fatou cylinder}. 
It parametrizes the orbits converging under the dynamics of the
iterate $f_c^m$ to the point $\zeta$. A repelling Fatou cylinder $C_R\simeq \CC/\ZZ$ is defined
similarly, as the quotient of a repelling petal.

Let $\tau$ be any conformal isomorphism $C_A\to C_R$. After uniformization,
$$C_A\underset{\approx}{\mapsto}\CC/\ZZ,\;C_R\underset{\approx}{\mapsto}\CC/\ZZ$$
$\tau(z)\equiv z+q\mod\ZZ$ for some $q\in\CC$.
Let $g_\tau:P_A\to P_R$ be any lift of $\tau$; it necessarily commutes with $f_c^{mq}$. 
Consider the semigroup $G$ generated by the dynamics of the pair
$(f_c,g_\tau)$. The orbit $Gz$ of a point $z\in\CC$ is independent of the choice
of the lift $g_\tau$ and only depends on $\tau$.

Set $$J_{(c,\tau)}=\{z\in\CC\text{ such that }Gz\cap J_c\neq\emptyset\}.$$
It can be shown that this set is the boundary of 
$$K_{(c,\tau)}=\{z\in\CC\text{ such that }Gz\text{ is bounded}\}.$$
Notice that $K_{(c,\tau)}\subsetneq K_c$: some of the orbits which converge to
$\zeta$ under $f_c$ are thrown into the complement $(\CC\setminus K_c)\cap P_R$ by $g_\tau$.
Holes which thus open in the set $K_c$ motivate the
use of the term ``implosion''.

\begin{figure}[ht]
\label{impl fig}
\centerline{\includegraphics[width=0.8\textwidth]{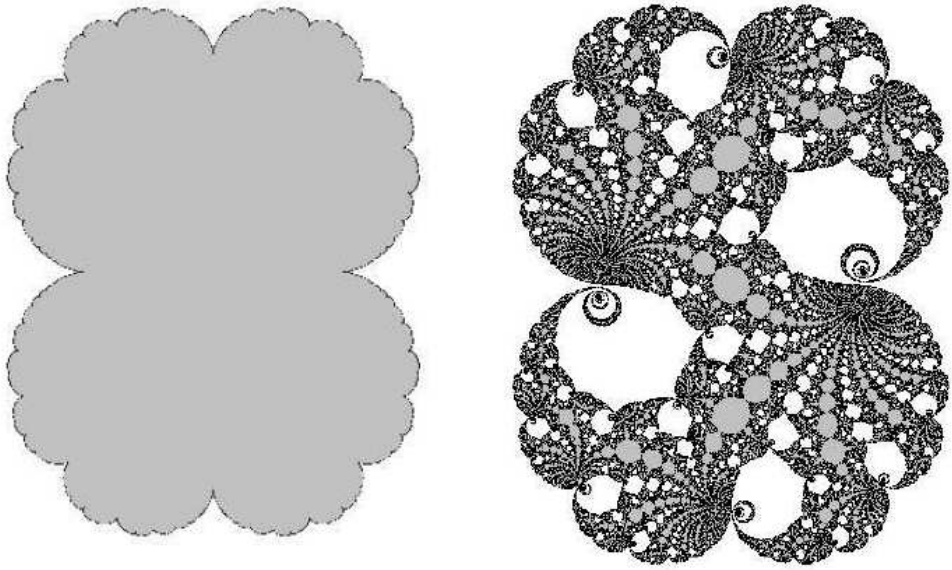}}
\caption{Before and after a parabolic implosion. The Julia sets (black) and filled Julia sets (light gray)
of a parabolic quadratic $f_{1/4}$ (left), and of $f_{1/4+\eps}$ for a small complex $\eps$.} 
\end{figure}

The Douady-Lavaurs theory postulates:

\begin{thm}
\label{implosion}
For every $\tau$ as above and every $s>0$ 
there exists $\tl c\in B(c,s)$ such that
$B(J_{\tl c},s)\supset J_{(c,\tau)}$.
\end{thm}

\noindent
Thus the Julia set of $f_c$ grows ``bigger'' under the perturbation from $c$ to $\tl c$.

\section{Positive results}
\label{sec:filled}

\subsection{Computability of filled Julia sets.}
In this section we show:
\begin{thm}
\label{thm:filled}
For any polynomial $p(z)$ there is an oracle Turing Machine
$M^{\phi}(n)$ that given an oracle access to the coefficients
of $p(z)$, outputs a $2^{-n}$-approximation of 
the filled Julia set $K_{p}\equiv K(p(z))$. 

\end{thm}

\noindent
Moreover,

\begin{thm}
\label{thm:filled quad}
In the case when $p(z)=z^2+c$ is quadratic, only two oracle
machines suffice to compute all non-parabolic filled Julia sets: one 
for $c\in \cM$, and one for $c\notin \cM$. 
\end{thm}

\noindent
Theorem \ref{thm:filled} answers in the affirmative the question posed to us by J.~Milnor,
after we first demonstrated the existence of non-computable quadratic Julia sets in 
\cite{BY}.

Let us first formulate a general fact:
\begin{prop}
\label{find roots}
Let $Q(z)$ be a complex polynomial. Then there exists a Turing Machine
$M^\phi$ with an oracle input for the coefficients of $Q(z)$ such that the
following holds. 
Consider any dyadic ball $B=B(\bar x,r)\subset \CC$, $\bar x\in\DD^2$, $r\in\DD$,
and  let $\alpha_1,\ldots,\alpha_m$
be the roots of $Q(z)$ contained in $B$. For any natural number $n$, the machine
$M^\phi$ will take $n$, $r$, and $\bar x$ as inputs, and will output a finite 
sequence of complex numbers $\beta_1,\ldots,\beta_k$ 
with dyadic rational real and imaginary parts
for which:
\begin{itemize}
\item $\beta_i\in B(\bar x,r+2^{-n} )$;
\item each $\beta_i$ lies at a distance not more than $2^{-n}$ from some root of $Q(z)$;
\item for every $\alpha_j$ there exists $\beta_i$ with $|\alpha_j-\beta_i|<2^{-n}$.

\end{itemize}

\end{prop}

\noindent
For a classical reference, see \cite{Wey}; a review of modern approaches to 
iterative root finding algorithms may be found in \cite{BCSS}.

For a given polynomial $p(z)$ we construct a machine computing 
the corresponding filled Julia set $K_p$. We will use 
some combinatorial information about $p$ in the construction, so the algorithm
will, in general, vary with the polynomial.
Note that all the information we will need can be encoded using a finite number 
of bits. 

\begin{itemize}
\item 
Information that would allow us to compute
the non-repelling orbits of the polynomial
with an arbitrary precision, as well as their type: attracting, parabolic, Siegel, or Cremer. 
By Fatou-Shishikura bound, there are at most $\deg p-1$ of them. 

By \propref{find roots}, such information could, for example, consist
of the list of periods $k_i$ of such orbits;
and for each $i$ a finite collection of dyadic
balls $\{D_i^j\}_{j=1}^{k_i}$ separating the points of the 
corresponding orbit from
the other solutions of the equation $p^{k_i}(z)=z$.

\item
For each (super)attracting periodic orbit 
$\bar\ze=\{\ze_1,\ldots,\ze_k\}$, a finite
union of dyadic balls $D_{\bar\ze}=\cup B(\ze_i,r_i)$ with the property
$$p^k(D_{\bar \ze})\Subset D_{\bar\ze}.$$

\item For each parabolic periodic point with period $m$ and multiplier $p/q$,
the values of $m$, $p$, $q$.

\item 
In the case of a Siegel disc $D$, information that would allow us
to identify a repelling periodic point $\ze_D$ in the same connected component
of $K_p$ as $D$. 
Again, by \propref{find roots}, it is sufficient to know its period, and a 
small enough dyadic ball around it, which separates it 
from all other points, periodic with the same period.
\end{itemize}

\subsection{Computing $K_p$}

We are given a dyadic point $d \in \DD$ and an $n\in \NN$. Our goal is 
to always terminate and output $1$ if $B(d,2^{-n}) \cap K_p \neq
\emptyset$ and to output $0$ if $B(d,2\cdot 2^{-n}) \cap K_p = \emptyset$. 
We do it by constructing five machines. They are guaranteed to 
terminate each on a different condition, always with a valid answer. 
Together they cover all the possible cases. 

\begin{lem}
\label{lem:filled}
There are five oracle machines $M_{ext}$, $M_{jul}$, $M_{attr}$, $M_{par}$, 
$M_{sieg}$ such that 
\begin{enumerate}
\item
if $d$ is at distance $\ge \frac{4}{3}\cdot 2^{-n}$ from $K_p$, 
$M_{ext}(d,n)$ will halt and output $0$. If $d$ is 
at distance $\le 2^{-n}$ from $K_p$, $M_{ext}(d,n)$ will 
never halt;
\item
if $d$ is at distance $\le \frac{5}{3} \cdot 2^{-n}$ from $J_p$, 
$M_{jul}(d,n)$ will halt and output $1$. If $d$ is at
distance $\ge 2\cdot 2^{-n}$ from $J_p$, $M_{jul}(d,n)$  will
never halt; 
\item 
$M_{attr}(d,n)$ halts and outputs $1$ if and only if $d$ is inside
the basin of an attracting orbit of $p$;
\item
$M_{par}(d,n)$ halts and outputs $1$ if and only if $d$ is inside
the basin of a parabolic orbit of $p$;
\item 
$M_{sieg}(d,n)$ halts and outputs
$1$ if the orbit of $d$ reaches a Siegel disc, and 
$d$ is at distance $\ge \frac{4}{3} \cdot 2^{-n}$ from $J_p$. 
It never halts if $d$ is at distance $\ge 2 \cdot 2^{-n}$
from $K_p$. 
\end{enumerate}
\end{lem}

\noindent

\begin{proof}[Proof of Theorem \ref{thm:filled}, given Lemma 
\ref{lem:filled}] By  Fatou-Sullivan classification 
 it is clear that for each $(d,n)$ at
least one of the machines halts. Moreover, by the definition 
of the machines, they always output a valid answer whenever they
halt. Hence running the machines in parallel and returning
the output of the first machine to halt gives the algorithm 
for computing $K_p$. 
\end{proof}

We now prove Lemma \ref{lem:filled}.

\begin{proof} {\bf (of Lemma \ref{lem:filled})}
We give a simple construction for each of the five machines. 
\begin{enumerate}
\item
$M_{ext}$: Take a large ball $B$ such that $p^{-1}(B) \Subset B$. 
Intuitively, we pull the ball back under $p$ to get a good 
approximation of $K_p$. Let $B_k$
be a $2^{-(n+3)}$-approximation of the set $p^{-k}(B)$. Output 
$0$ iff $B_k \cap B(d,\frac{7}{6}\cdot 2^{-n}) = \emptyset$. 
It is not hard to see that this algorithm satisfies the conditions 
on $M_{ext}$. 
\item 
$M_{jul}$: 
By \propref{find roots} for each $k$ we can compute all periodic orbits
of $p(z)$ in $B$ with periods $j\leq k$, as roots of the equation
$$p^j(z)-z=0$$
with an arbitrarily high precision.
Moreover, by our assumptions, we have the means to distinguish the 
non-repelling orbits from the repelling ones.

Let $C_k$ be a finite collection of complex numbers with dyadic rational
real and imaginary parts which approximate the  repelling periodic orbits  
with periods up to $k$ with precision $2^{-(n+3)}$.
Output $1$ iff $d(d,C_k) < \frac{11}{6}\cdot 2^{-n}$.
The repelling periodic orbits are all in $J_p$ and are dense in
this set. Hence the algorithm satisfies the conditions on $M_{jul}$. 
\item 
$M_{attr}$: For each attracting orbit $\bar\zeta$ of period $m_{\bar\zeta}$
 find $l_{\bar\zeta}$ such that
$$B(p^{m_{\bar\zeta}}(D_{\bar\zeta}),{2^{-l_{\bar\zeta}}})\subset D_{\bar\zeta}.$$
Let 
$$l=1+\sup l_{\bar\zeta},\text{ and }m=\prod m_{\bar\zeta}.$$ 
Let $z_k$ be a $2^{-(l+3)}$-approximation of $p^{mk}(d)$. 
If $d$ is inside the  basin of  an  attracting orbit $\bar\ze$, then 
$z_k$ will be inside $D_{\bar\ze}$ for some $k$. Output $1$ 
if $z_k$ is inside $D_{\bar\zeta}$ and at least $2^{-l}$-far from the boundary of $D_{\bar\ze}$. 
\item 
$M_{par}$: 
We make use of \lemref{constructing petal}. Since we can produce arbitrarily good approximations
of every parabolic periodic point $\zeta$ of $p(z)$, we do not need an oracle for the value of this
point. Let $L_k^\zeta$ be the sets from \lemref{constructing petal} corresponding to the point 
$\zeta$. Let $z_k=p^k(d)$ computed with precision $2^{-(k+2)}$. We output $1$ if $z_k$ is inside
$L_k^\zeta$ for some $\zeta$ and at least $2^{-k}$-away from its boundary.
\item 
$M_{sieg}$: This is the most interesting case. It is not 
hard to see that for each $k$, we can compute a union $E_k$ of 
dyadic balls such that 
$$
\bigcup_{i=0}^k p^i \left(B(d,\frac{4}{3} \cdot 2^{-n})\right) \subset 
E_k \subset \bigcup_{i=0}^k p^i \left(B(d, \frac{5}{3} \cdot 
2^{-n})\right).
$$
Let $\zeta_*$ be the center of the Siegel disc (one of the centers, in case 
of an orbit), and let $y$ be the given periodic point in the connected 
component of $\zeta_*$. We terminate and output $1$ if $E_k$ separates $\zeta_*$ from 
$y$ in $\CC$ (or covers either one of them) for some $k$. 

Clearly, if $d$ is inside the Siegel disc, then the forward images of
$B(d, \frac{4}{3} \cdot 2^{-n})$ will cover an annulus in the 
disc that will separate $\zeta_*$ from the boundary of the disc, and
in particular from $y$. Hence $M_{sieg}$ will terminate and output 
$1$. 

On the other hand, if the distance from $d$ to $K_p$ is $\ge 2\cdot 2^{-n}$, then
$E_k \cap K_p = \emptyset$ for all $k$. In particular, $E_k$ 
cannot separate $\zeta_*$ from $y$, since they are connected in $K_p$.
\end{enumerate}
\end{proof}

\noindent
The proof is simplified in the case of a quadratic polynomial.

\begin{proof}[Proof of \thmref{thm:filled quad}].
If we assume that $p(z)=f_c(z)$ then by the Fatou-Shishikura bound,
there is at most one non-repelling orbit. By our assumption, it is not parabolic.
Moreover, if it is a Siegel orbit, then the Julia set is connected. Therefore,
{\it any} repelling periodic orbit will be in the same connected component
of $K_p$ as the Siegel disk.

If $c\notin \cM$, we run $M_{ext}$ and $M_{jul}$. One and only one of them is
guaranteed to halt and output a correct answer.

For $c\in\cM$ we will use a modified Turing Machine $\widehat M_{sieg}.$ It will
compute the set $E_k$ as before. 
If $E_k$ separates the plane, it will use \propref{find roots} to search 
for a periodic point of period at most $k$ both in the exterior and the
interior components of $E_k$. If a $k$ is found for which such two orbits are located,
or if $E_k$ covers a periodic orbit,
it will terminate and output $1$.

If $c\in \cM$, then we run $M_{ext}$, $M_{jul}$, $M_{hyp}$, and $\widehat M_{sieg}$. As before,
it is easy to see that one of them will terminate, and its output will be a correct one.
\end{proof}

\begin{cor}
Denote by $\cP$ the set of 
$c$'s for which $J_c$ is parabolic. 
The function $K: c \mapsto K_{z^2+c}$ is continuous in the 
Hausdorff metric on the set $\cM\setminus\cP$. 
\end{cor}

\subsection{Computability of Julia sets in the absense of rotation domains}
Similar ideas were used in \cite{BBY1} to prove the following theorem:

\begin{thm}
\label{comp rat jul}
Let $f$ be a rational map $f:\hat \CC\to\hat\CC$ without rotation domains.
Then its Julia set is computable in the spherical metric by an oracle Turing machine $M^\phi$ with 
the oracle representing the coefficients of $f$. The algorithm uses the following non-uniform information
about each parabolic periodic point $\zeta$ of $f$ with period $m$ and multiplier $e^{2\pi i p/q}$:
\begin{itemize}
\item a dyadic ball $B(w,r)\ni p$ such that $B(w,2r)$ does not contain any other points periodic with period
$m$;
\item the values of $m$, $p$, and $q$.
\end{itemize}
\end{thm}

\begin{proof}
For every natural $n$ we can compute a sequence of rationals $\{q_i\}$ 
such that 
\beq
\label{lbeq}
B(J_f,2^{-(n+2)})  \Subset \bigcup_{i=1}^{\infty} B(q_i, 2^{-(n+1)}) 
\Subset B(J_f, 2^{-n}).
\eeq

To do that, for each $k>n+2$ we compute $2^{-k}$-approximations of the periodic points
of $f$ in $\hat\CC$ with periods at most $k$ using \propref{find roots}.
Let $M>0$ be some bound on $|D^2 f^m(z)|$ in the area of an approximate periodic orbit  $r_i$ with period $m$. 
Then $|Df^m(r_i)|>1+2^{-k}M$ 
implies that $|Df^m(w)|>1$
for the periodic point $w$ which $r_i$ approximates. In this case we add the point $r_i$ to our sequence
of rationals. Clearly, for each repelling periodic point of $f$ we will eventually obtain in this way
a rational point which approximates it with precision at least $2^{-(n+3)}$. Since 
such points are contained in $J_f$, and dense there, our sequence has the desired property.

Of course, we can similarly {\it eventually} find every attracting orbit $\bar\zeta$ of $f$ with an arbitrary
precision. In this case, we will compute a set $D_{\bar\zeta}$ for this orbit with the same properties as before.
Set $D=\cup_{\bar\zeta}D_{\bar\zeta}$.

Finally, for each parabolic periodic point $\zeta$ of $f$ let $L_k^\zeta$ be the sets from \lemref{constructing petal}.
Set $L_k=\cup_\zeta L_k^\zeta$.

We are now ready to present an algorithm to find a set $C_m\in\cC$ with $\dist_H(C_m,J_f)<2^{-m}$.
Fix $m\in\NN$. Our algorithm to find $C_m\in\cC$ works as follows.
At the $k$-th step:

\begin{itemize}
\item  compute the finite union $B_k=\cup_{i=1}^k B(q_i,2^{-(m+1)})\in \cC$;
\item compute with precision $2^{-(m+3)}$ the complement of the preimage 
$$f^{-k}(D\cup L_k),$$ that is,
find $W_k\in \cC$ such that 
$$d_H(\overline{W_k},\overline{\hat\CC\setminus (f^{-k}(D\cup L_k)}))<2^{-(m+3)};$$
\item if $W_k\subset B_k$ output $C_m=B_k$ and terminate. Otherwise, go to step $k+1$.
\end{itemize}

\noindent
By Fatou-Sullivan classification, the algorithm will eventually terminate. 
Now suppose that the algorithm terminates on step $k$. Since $W_k\subset B_k$
and $J_f\subset B(W_k,2^{-(m+3)})$ we have $J_f\subset B(C_m,2^{-(m+3)})$. On the other hand,
$\cup \{ q_i\}\subset J_f$, and thus $B_k=C_m\subset B(J_f,2^{-(m+1)})$.

\end{proof}

\section{Computability of Julia sets of Siegel quadratics and negative results}
\label{sec:siegel}
\subsection{Computabilty of $r(\theta)$ is equivalent to computability of $J_\theta$}.
Recall the discussion of the family $P_\theta(z)=z^2+e^{2\pi i\theta}z$ from \secref{sec:quad siegel}.
As before, if $\theta$ is a Brjuno number, we denote $\Delta_\theta$ its Siegel disk,
and $r(\theta)$ the conformal radius of $\Delta_\theta$. When $\theta\notin\cB$, we set 
$r(\theta)=0$. Let us denote $r_{\text{sup}}=\sup r(\theta)$. We let $J_\theta$ stand for $J(P_\theta)$.

Of course, the change in parametrization from $c$ to $\theta$ makes it natural to talk
about computabilty of $J_\theta$ by a TM with an oracle for $\theta$, rather than for $c$.
However, these notions are obviously equivalent, as $c=c(\theta)$ is found by the formula:

\beq
\label{comp c theta}
c=c(\theta)=\lambda/2-\lambda^2/4,\text{ where }\lambda=e^{2\pi i\theta}.
\eeq

To address the question of computability of $J_\theta$ for $\theta\in\cB$ we first make note
of the following result, proven in \cite{BBY1}:

\begin{prop}
\label{prop:comp2}
Suppose $r_\theta$ is computable by a Turing Machine $M^\phi$ with an oracle for $\theta$.
Then so is $J_{\theta}$.
\end{prop}

To clarify the logic of the argument, let us break the proof into two steps:
\begin{lem}
Suppose $r_\theta$ is computable by a Turing Machine $M^\phi$ with an oracle for $\theta$.
Then so is the inner raidus $\rho(\Delta_\theta,0)\equiv \rho_\theta$.
\end{lem}
\begin{proof}
The algorithm works as 
follows:
\begin{itemize}
\item[(I)] For $k\in\NN$  compute a set $D_k\in\cC$ which is a 
$2^{-m}$-approximation of the preimage $P_\theta^{-k}(D)$, for some sufficiently
large disk $D$;
\item[(II)]  evaluate the
conformal radius $r(D_k,0)$ with precision $2^{-(m+1)}$ (this can be done, for example,
by using one of the numerous existing methods for computing the Riemann Mapping of a computable domain,
see \cite{BB});
\item[(III)] as before, denote $$F(x)=4x/(1+x)^2\text{, for }x\in[0,1].$$
 Note that this function is monotone, and
let $\psi(w)=F^{-1}(w)$. This function is computable, and $\psi(1)=1$.

Evaluate 
$$p=\psi(r_\theta/r(D_k,0))$$
 with precision $2^{-(m+5)}/\rho(D,0)$. If 
$$|1-p|<2^{-(m+3)}/\rho(D,0),$$
then compute the inner radius
$\rho(D_k,0)\equiv r_k$ around $0$ with precision $2^{-{(m+1)}}$ 
 and output this number. Else, increment $k$ and
return to step (I).

\end{itemize}

\medskip
\noindent
{\bf Termination.}
Let $K=K({P_\theta})$ be the filled Julia set of $P_\theta$. 
Then $$\cap_{k=0}^{\infty} D_k = K \supset \Delta_\theta$$
and $D_0 \supset D_1 \supset D_2 \supset \ldots$. 

Hence for every $\delta>0$ there will be a step $k=k(\eps)$ after which
$$\dist(\partial D_k,J_{\theta})<\delta.$$
Since $\partial \Delta_\theta \subset J_{\theta}$,  by \lemref{variation conf radius} this implies that 
$$|r(D_k,0)-r(\Delta_\theta,0)|=|r(D_k,0)-r_\theta|<4\sqrt{r(D,0)}\sqrt{\delta}\underset{\delta\to 0}{\longrightarrow} 0.$$
Since for every large enough  $k$, the value of 
$$\psi(r_\theta/r(D_k,0))>1-2^{-(m+4)}/\rho(D,0),$$
the algorithm will eventually terminate on step (III).

\medskip
\noindent
{\bf Correctness.}
Now suppose the algorithm has terminated on step (III). 
As $\Delta_\theta \subset D_k$,  \lemref{variation conf radius} implies that 
$$1\geq \frac{\rho_\theta}{\rho(D_k,0)}\geq 1-\frac{2^{-(m+1)}}{\rho(D,0)},$$
and so
$$|\rho(D_k,0)-\rho_\theta|\leq 2^{-(m+1)}.$$

\end{proof}

\begin{lem}
\label{rho comput}
Suppose $\rho_\theta$ is computable by a Turing Machine $M^\phi$ with an oracle for $\theta$.
Then so is $J_{\theta}$.
\end{lem}
\begin{pf}
The algorithm to produce the $2^{-n}$ approximation of the Julia set is the following.
First, compute a large disk $D$ around $0$ with $P_\theta(D)\Supset D$. Then,
\begin{itemize}
\item[(I)] compute a set $D_k\in\cC$ which is a 
$2^{-(n+3)}$-approximation of the preimage $P_\theta^{-k}(D)$;
\item[(II)] set $W_k$ to be the round disk with radius $\rho_\theta-2^{-k}$ about the
origin. Compute a set $B_k\in \cC$ which 
is a  $2^{-(n+3)}$-approximation of $P_\theta^{-k}(W_k)$;
\item[(III)] if $D_k$ is contained in a $2^{-(n+1)}$-neighborhood of $B_k$, then output
a $2^{-(n+1)}$-neighborhood of $D_k \setminus B_k$, and stop. If not, go to step (I).
\end{itemize} 

\noindent
A proof of the validity of the algorithm is obvious, and we leave it to the reader.

\end{pf}

\begin{figure}[ht]
\centerline{\label{fig-siegel2}
\includegraphics[width=6cm]{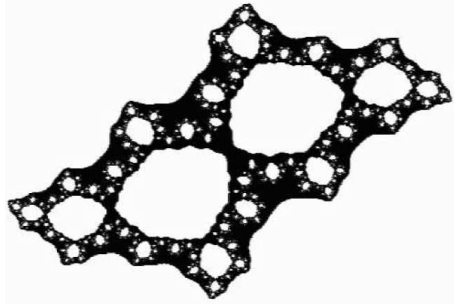}}
\caption{A figure produced by the algorithm of \lemref{rho comput} for $\theta=(\sqrt{5}-1)/2$.
It does not have the same artistic quality as \figref{fig-siegel}, but has a guaranteed accuracy
up to the selected size of a pixel.}
\end{figure}

\noindent
By \propref{prop:comp2} and \corref{noncomp radius} we have:
\begin{thm}
\label{rad computable}
The conformal radius $r(\theta)$ is computable by a Turing Machine with an oracle for $\theta$ if and only
the same is true for the Julia set $J_\theta$.

\end{thm}

\noindent
Let us make a note:
\begin{prop}
\label{bded comp}
Let $\theta$ be of bounded type. Then $J_\theta$ is computable by a TM with an oracle for $\theta$.
\end{prop}
\begin{pf}
By \propref{bdary bded type}, $\rho_\theta$ is computable. The claim follows by \lemref{rho comput}.
\end{pf}

\noindent
Showing that there exist {\it noncomputable} Siegel Julia sets is much more delicate. In the next
section we will prove that this can happen even if $\theta$ itself is computable.

\subsection{Conformal radius of a Siegel quadratic with a computable $\theta$}
The theorem we formulate below characterizes the values of $r(\theta)$ which correspond
to {\it computable} parameters $\theta$:

\begin{thm}
\label{thm:char} 
Let $r \in (0,r_{\text{sup}})$ be a real number. Then $r=r(\theta)$ is the conformal 
radius of a Siegel disc of the Julia set $J_\theta$ for some computable 
number $\theta$ if and only if $r$ is right-computable.
\end{thm}

\noindent
Before proving this theorem, let us formulate a corollary:

\begin{cor}
\label{cor:sieg}
There exist computable values of parameter $c$, such that the Julia set $J_c$ is not computable
by a TM $M^\phi$ with an oracle access to $c$.
\end{cor}

\begin{proof}
By \propref{right-comp} there exist  right computable numbers $r_*\in[0,r_{\text{sup}}]$
which are not in $\RR_\cC$. By \thmref{thm:char}, $r_*=r(\theta_*)$ for $\theta_*\in\RR_\cC$.
Since $\theta_*$ itself is computable, $r_*$ is uncomputable by a TM with an oracle access 
to $\theta_*$. By \thmref{rad computable}, the Julia set $J_{\theta_*}$ is uncomputable by
a TM with an oracle access to $\theta_*$. The claim follows by 
(\ref{comp c theta}).
\end{proof}

\begin{proof}[Proof of the ``only if'' direction of \thmref{thm:char}]
We assume that $\theta$ is computable, and show that $r(\theta)$ is right-computable. 
Recall that periodic orbits are dense in the Julia set $J_\theta$. Let $H_n$ 
be the union of all repelling periodic orbits with periods between $1$ and $n$.
By \propref{find roots} we can algorithmically find an arbitrarily good approximation
of $H_n$ by dyadic rationals.

$J_\theta$ is connected, and $\cup H_n$ is dense in $J_\theta$. 
Thus for every $l$, there exists $n_l$ such that 
the set $B(H_{j},2^{-(l+1)})$ is connected. Moreover, such $n_l$ can 
be found algorithmically. 

Since $J_\theta$ separates $\al$ from $\infty$, the same is true for 
$B(H_{n_l},2^{-(l+1)})$ provided $l$ is sufficiently large.
Hence we can compute a strictly increasing sequence $\{n_l\}_{l=l_0}^{\infty}\subset\NN$,
and a set $U_l\subset \cC$ with the property
$$
B(H_{n_l},2^{-(l+1)})\subset U_l\subset B(H_{n_l},2^{-l}).
$$
such that $\CC\setminus \bar U_l$ has a simply-connected component $W_l$ containing $\al$.

 Using any constructive algorithm for computing the conformal radius \cite{BB}
we can approximate the $k$-th term of the sequence 
$$
R_k = r(B(W_k,2^{-(k-1)}),\al)+5\cdot 2^{4-\frac{k-1}{2}}
$$
$H_n \ra J_\theta$ in Hausdorff metric and $n_k \ra \infty$, thus by \lemref{variation conf radius},  $R_k \ra r(\theta)$. Moreover, 
$\{R_k\}$ is a non-increasing sequence. Let 
$\rho_k$ be a dyadic approximation of $R_k$ that we compute so that $|\rho_k - R_k|<2^{-k}$. 
Let 
$$
r_k = \rho_k + 3 \cdot 2^{-k}. 
$$
Then $\{r_k\}$ is a computable sequence of dyadic numbers. We have
$$
\lim_{k\ra \infty} r_k = \lim _{k\ra \infty} \rho_k = \lim _{k\ra \infty} R_k = r(\theta),
$$
and for each $k$,
$$
r_k = \rho_k + 3 \cdot 2^{-k} \ge R_k + 2\cdot 2^{-k} \ge R_{k+1} + 4 \cdot 2^{-(k+1)} \ge
\rho_{k+1} + 3\cdot 2^{-(k+1)} = r_{k+1}.
$$
This shows that $r(\theta)$ is right-computable.

 Note that we know that $H_n\ra J_\theta$ is 
Hausdorff metric, which allows us to conclude that $r_k\ra r(\theta)$. However, we do not have 
(and {\em cannot} have) an estimate on the rate of convergence of $H_n$ to $J_\theta$, and 
thus cannot obtain an estimate on the rate of convergence of $r_k\ra r(\theta)$ and compute $r(\theta)$. 
\end{proof}

\subsection*{Proof of the ``if'' direction of \thmref{thm:char}}
Given a computable sequence $\{r_n\}$ such that $r_n \searrow r$ we 
claim that we can construct a $\theta$ such that $r= r(\theta)$.  We will be using 
the following three lemmas. The first one is Lemma~3.1 of \cite{BY}, and the second
one is Lemma~4.2 of \cite{BBY2}. The proofs are outlined in Appendix \ref{app:B}.

\begin{lem}
\label{smlchg}
For any initial segment $I = [a_0, a_1, \ldots, a_n]$, write 
$\omega = [a_0, a_1, \ldots, a_n, 1, 1, 1, \dots]$. Then 
for any $\ve>0$, there is an $m>0$ and an integer $N$
such that if we write $\be = [a_0, a_1, \ldots, a_n, 1, 1, \ldots, 
1, N, 1, 1, \ldots]$, where the $N$ is located in the $n+m$-th 
position, then
$$
\Phi(\omega) + \ve < \Phi(\be) < \Phi(\omega) + 2 \ve. 
$$
\end{lem}

\begin{lem}
\label{notdeclem}
For  $\omega$ as above, for any $\ve>0$ there is an $m_0>0$, which 
can be computed from $(a_0, a_1, \ldots, a_n)$ and $\ve$,  
such that for any $m \ge m_0$, and for any tail 
$I = [a_{n+m}, a_{n+m+1}, \ldots]$ if we denote 
$$ \be^{I} = [a_1, a_2, \ldots, a_n, 1, 1, \ldots, 1, a_{n+m}, a_{n+m+1}, 
\ldots],$$ then
$$
\Phi(\be^{I}) > \Phi (\om) - \ve. 
$$
\end{lem}

\begin{lem}
\label{lem:tailof1s}
Let $\om = [a_1, a_2, \ldots]$ be a Brjuno number, that is $\Phi(\om)<\infty$. 
Denote $\om_k = [a_1, a_2, \ldots, a_k, 1,1, \ldots]$. 
Then for every $\ve>0$ there is an $m$ such that for all $k\ge m$,
$$
\Phi(\om_k) < \Phi(\om)+\ve. 
$$
\end{lem}

Using \propref{bded comp}, we can get a computable version of Lemmas \ref{smlchg} 
and \ref{notdeclem}. 

\begin{lem}
\label{smlchgr}
For any given initial segment $I = [a_0, a_1, \ldots, a_n]$ and $m_0>0$, write 
$\omega = [a_0, a_1, \ldots, a_n, 1, 1, 1, \dots]$. Then 
for any $\ve>0$, we can uniformly compute  $m>m_0$, an integer $t$ and an integer $N$
such that if we write $\be = [a_0, a_1, \ldots, a_n, 1, 1, \ldots, 
1, N, 1, 1, \ldots]$, where the $N$ is located in the $n+m$-th 
position, we have 
\beq
\label{rcond}
r(\omega) - 2 \ve < r(\be) < r(\omega) -  \ve,
\eeq
\beq
\label{phiinc}
\Phi(\be)> \Phi(\om),
\eeq
and for any
$$
\ga = [a_0, a_1, \ldots, a_n, 1, 1, \ldots, 1, N, 1, \ldots, 1, c_{n+m+t+1}, c_{n+m+t+2}, \ldots],
$$
\beq
\Phi(\ga)>\Phi(\om)-2^{-n}. 
\eeq
\end{lem}

\begin{proof}
We first show that such $m$ and $N$ {\em exist}, and then 
give an algorithm to compute them. By Lemma \ref{smlchg}
we can increase $\Phi(\omega)$ by any controlled amount 
by modifying one term arbitrarily far in the expansion. 

By Theorem \ref{phi-cont}, $f: \theta\mapsto \Phi(\theta)+\log r(\theta)$
extends to a continuous function. Hence for any $\ve_0$ there is a $\de$ such 
that $|f(x)-f(y)|<\ve_0$ whenever $|x-y|<\de$. In particular, 
there is an $m_1$ such that $|f(\be)-f(\om)| < \ve_0$ whenever $m \ge m_1$.

This means that if we choose $m$ large enough, a controlled increase 
of $\Phi$ closely corresponds to a controlled drop of $r$ by a corresponding 
amount, hence there are $m>m_0$ and $N$ such that 
\eref{rcond} holds. \eref{phiinc} is satisfied almost automatically.
The only problem is to {\em computably} find such $m$ and $N$. 

To this end, we apply \propref{bded comp}. Together with \thmref{rad computable},
it implies that for any specific $m$ 
and $N$ we can compute $r(\be)$. This means that we can find the suitable 
$m$ and $N$, by enumerating all the pairs $(m,N)$ and exhaustively checking 
\eref{rcond} and \eref{phiinc} for all of them. We know that eventually we will find a 
pair for which \eref{rcond} and \eref{phiinc} hold. 

Finally, $t$ exists and can be computed by Lemma \ref{notdeclem}.
\end{proof}

Lemma \ref{lem:tailof1s} yields the following lemma. 

\begin{lem}
\label{lem:supsup}
The supremum of $r(\theta)$ over all angles is equal to the supremum 
over the angles whose continued fraction expansion has only finitely 
many terms that are not $1$:
$$
r_{sup} = \sup_{\theta=[a_1,a_2,\ldots,a_k,1,1,\ldots]} r(\theta).
$$
\end{lem}

\begin{proof}
Let $\ve>0$ be an arbitrary small positive number. By the definition 
of $r_{sup}$ there is a $\theta = [a_1,a_2, \ldots]$ such that 
$\log r(\theta)>\log r_{sup}-\ve$. Denote
$$\theta_k = [a_1, a_2, \ldots, a_k, 1,1,\ldots].$$
Lemma \ref{lem:tailof1s} states that there is an $m$ such that for $k\ge m$, $\Phi(\theta_k)<
\Phi(\theta)+\ve$. Moreover, there is a $\de$ such that whenever 
$|\phi-\theta|<\de$ we have $|\upsilon(\phi)-\upsilon(\theta)|< \ve$. 

$\theta_k \ra \theta$, hence there is an $n\ge m$ such that $|\theta_n - \theta|<\de$. 
$\theta_n$ has the required form, and we have
$$
\log r(\theta_n) = \upsilon(\theta_n) - \Phi(\theta_n) > \upsilon(\theta)-\Phi(\theta)-2 \ve = 
\log r(\theta) - 2 \ve > \log r_{sup} - 3 \ve. 
$$
This shows that we can make $r(\theta_n)$ as close to $r_{sup}$ as we like. 
\end{proof}

We are given $r = \lim \searrow r_n < r_{sup}$, hence there is an $s$ and an $\ve>0$ such
that $r_s < r_{sup} - 2 \ve$. By Lemma \ref{lem:supsup}, there is
a $\ga_0=[a_1,a_2,\ldots,a_n,1,1,\ldots]$ such that $r_s+\ve/2 < r(\ga_0)<r_{s}+\ve$. 

We are now ready to give an algorithm for computing a rotation number $\theta$ 
for which $r(\theta)=\lim \searrow r_n$. 
The algorithm works as follows. On stage $k$ it produces a finite initial segment 
$I_k=[a_0, \ldots, a_{m_k}]$ such that the following properties are maintained:
\begin{enumerate}
\item 
$I_0=[a_1,a_2,\ldots,a_n]$;
\item 
$I_k$ has at least $k$ terms, i.e. $m_k \ge k$;
\item 
for each $k$, $I_{k+1}$ is an extension of $I_k$;
\item 
for each $k$, denote $\ga_k = [I_k, 1, 1, \ldots]$, then $r_{s+k}+2^{-(k+1)} \ve< r(\ga_k) < r_{s+k} + 2^{-k} \ve$; 
\item
for each $k$, $\Phi(\ga_k)> \Phi(\ga_{k-1})$;
\item 
for each $k$, for any extension 
$$
\be = [I_k, b_{m_k+1}, b_{m_k+2}, \ldots], 
$$
$\Phi(\be) > \Phi(\ga_k)-2^{-k}$. 
\end{enumerate}

The first three properties are very easy to assure. The last three are maintained 
using Lemma \ref{smlchgr}. By this Lemma we can decrease $r(\ga_{k-1})$ by any given 
amount (possibly in more than one step) by extending $I_{k-1}$ to $I_k$. Here we 
use the facts that the $r_k$'s are computable and non-increasing. 

Denote $\theta = \lim_{k\ra \infty} \ga_k$. The continued fraction expansion of 
$\theta$ is the limit of the initial segments $I_k$. 
This algorithm gives us at least one term of the continued fraction expansion 
of $\theta$ per iteration, hence we would need at most $O(n)$ iterations 
to compute $\theta$ with precision $2^{-n}$ (in fact, much fewer iterations 
would suffice). The initial segment $\ga_0$ can also be computed as in the proof 
of Lemma \ref{smlchgr}.
It remains to prove that, in fact, $\theta$ is the rotation number we are looking 
for. 

\begin{lem}
The following equalities hold:
$$
\Phi(\theta) = \lim_{k \rightarrow \infty}\Phi(\gamma_k)~~~~\mbox{$~$and$~$}~~~~
r(\theta) = \lim_{k \rightarrow \infty}r(\gamma_k)=r.
$$
\end{lem}

\begin{proof}
By the construction, the limit $\theta = \lim \ga_k$ exists. We also 
know that the sequence $r(\ga_k)$ converges  to 
the number $r=\lim\searrow r_k$, and that the sequence $\Phi(\ga_k)$ is 
monotone non-decreasing, and hence converges to a value $\psi$ ({\em a priori}
we could have $\psi=\infty$). By the Carath{\'e}odory Kernel Theorem, we have $r(\theta) \ge r>0$, so 
$\Phi(\theta) < \infty$.
On the other hand, by the property we have maintained through 
the construction, we know that $\Phi(\theta)>\Phi(\ga_k)-2^{-k}$ for all $k$. 
Hence $\Phi(\theta) \ge \psi$.  In particular, $\psi < \infty$. 

From \cite{BC} we know that 
\beq
\label{BCCont}
\psi + \log r= \lim (\Phi(\ga_k) + \log r(\ga_k)) = \Phi(\theta) + \log r(\theta).
\eeq
Hence we must have $\psi = \Phi (\theta)$, and $r = r(\theta)$, which completes
the proof. 
\end{proof}

\section{Interpretation of the results}
\label{sec:interpret}

\subsection{How difficult is it to produce a $\theta$ for which $J_\theta$ is uncomputable?}
As we have seen, a value of $\theta$ for which $J_\theta$ is uncomputable can be produced
constructively. As we will see below,  under a reasonable 
assumption, it is not even hard to do so:

\begin{tet}
\label{thm:poly}Assume that
the 1-periodic continuous function $\upsilon: \theta \mapsto \Phi(\theta) + \log r(\theta)$ has 
a computable modulus of continuity (\ref{mod cont}); this follows, for instance, from Marmi-Moussa-Yoccoz Conjecture. 
Suppose there is a computable sequence $r_1, r_2, \ldots$ of dyadic numbers such 
that 
\begin{itemize}
\item 
$\{r_i\}$ is non-increasing, $r_1 \ge r_2 \ge \ldots$, and 
\item 
$\lim_{i\ra \infty} r_i = r$. 
\end{itemize}
Then there is a {\em poly-time} computable $\theta$ (and hence a poly-time computable $c=c(\theta)$)
such that $r(\theta) = r$.
\end{tet}

\begin{proof}
By the assumption,  there is a computable function 
$\mu : \NN \ra \NN$ such that 
$$| \upsilon(\theta_1)-\upsilon(\theta_2)| < 2^{-n}\text{ whenever }| \theta_1 - \theta_2|<2^{-\mu(n)}.$$

The proof goes along the lines of the proof of the ``if" direction of Theorem \ref{thm:char}.
We outline the modifications made to the proof here and leave the details to the reader. 
The key difference is that in the proof of Theorem \ref{thm:char} we used Lemma \ref{smlchgr}
to perform a step in decreasing the conformal radius from $r(\ga_{k-1})$ to $r(\ga_k)$. 
The algorithm there is basically an exhaustive search, which, of course, could take much
longer than polynomial time in the precision of $\ga_k$ to compute. By assuming 
that $\upsilon$ has a computable modulus of continuity,
 we can deal with $\Phi(\ga_{k-1})$ and $\Phi(\ga_k)$ instead 
of the $r(\bullet)$'s. We have an explicit formula for $\Phi$ that converges well, and 
we can compute the continued fractions coefficients to make $\Phi(\ga_k)$ close to 
whatever we want relatively fast.

The step of going from $\ga_{k-1}$ to $\ga_k$ is as follows. First, we 
do the following computations:
\begin{itemize}
\item 
compute $d_k$ which is the ``drop" in $r$ we are trying to achieve; 
we want $$d_k/2 < \log(r(\ga_{k-1})) - \log(r(\ga_k)) < d_k;$$
\item 
compute using the function $\mu$ a value $\de_k$ such that $|\upsilon(x)-\upsilon(y)|<d_k/8$
whenever $|x-y|<\de_k$. 
\end{itemize}
We have no {\em a priori} bound on how long these computations would take, but 
we would still like to be computing $\theta$ in polynomial time. 
To achieve this, we use $1$'s in the continued fraction expansion of 
$\theta$ to ``pad" the computation. 

When asked about the value of $\theta$ with precision $2^{-n}$ which is higher 
than what the known terms of the expansion $[I_{k-1}]$ can provide, we do the following:
\begin{itemize}
\item 
try to compute $d_k$ and $\de_k$ as above, but run the computation for 
at most $n$ steps;
\item
if the computation does not terminate, output an answer consistent with 
the initial segment $[I_{k-1},\underbrace{1,1,\ldots,1}_{2 n}]$;
\item 
if the computation terminates in less than $n$ steps proceed as described 
below.
\end{itemize}
Note that so far the computation is polynomial in $n$. 
For some sufficiently large $n$ the computation will terminate in $n$ steps, 
at which point we will have computed $d_k$ and $\de_k$. 
If necessary,
we then add more $1$'s to the 
initial segment to assure that $|\ga_{k-1}-\ga_{k}|<\de_k$. 

Recall that our goal is to assure that 
$$d_k/2 < \log(r(\ga_{k-1})) - \log(r(\ga_k)) < d_k.$$
With the current initial segment for $\ga_k$ we have $|\ga_{k-1}-\ga_{k}|<\de_k$, and 
hence in the difference
$$
\log(r(\ga_{k-1})) - \log(r(\ga_k))  = 
\Phi(\ga_k) - \Phi(\ga_{k-1}) + (\upsilon(\ga_{k-1}) - \upsilon(\ga_k))
$$
the last term is bounded by $d_k/8$. 
This means that for the current step it suffices to increase $\Phi(\ga_k)$ relative to 
$\Phi(\ga_{k-1})$ by between $\frac{5}{8} d_k$ and $\frac{7}{8} d_k$. 

Let $M$ be the total length of $I_{k-1}$ and the $1$'s we have added, and let us extend the
continued fraction by putting $N\in\NN$ in the $M+1$-st term, and all $1$'s further.
Increasing $M$ if necessary, we can ensure an approximate equality 
$$
\Phi(\ga_k) \approx \Phi(\ga_{k-1}) + \al(N) \log N
$$
up to an error of $\frac{1}{32}d_k$.
Let $p_M/q_M$ be 
the $M$-th convergent of the resulting continued fraction. 
Recall that on an input $n$ we need to compute $\theta$ with 
precision $2^{-n}$ in time polynomial in $n$. If $2^{-n}> 1/\sqrt{q_M}$, 
then we do not need to know anything about $N$ to compute the required 
approximation. Suppose $2^{-n}<1/\sqrt{q_M}$, which means $n> \log q_M/2$. 
And we have time polynomial in $\log q_M$ to perform the computation. 

Note that $M = O(\log q_M)$. It is also not hard to see that $\al(N)<2^{-M/2}$,
so in order to have a change by $\approx 3 d_k/4$ we must have $N>e^{\Om(2^{M/2})}$, 
hence by making $M$ sufficiently large (depending on the value of $d_k$), 
we can guarantee that $N>e^{2^{M/3}}$. This means that we can 
approximate $\al(N)$ with the truncated function $\Phi$ at the $M$-th
convergent of the continued fraction. Write $p_M/q_M = [a_1, a_2, \ldots, a_M]$, and denote
$$
\be = [a_1, a_2, \ldots, a_M] \cdot [a_2, a_3, \ldots, a_M] \cdot \ldots
\cdot [a_{M-1}, a_M] \cdot [a_M].
$$
Then $\be$ approximates $\al(N)$ within a very small relative error. In particular, we can 
assure that 
$$
\be \cdot \left( 1- \frac{1}{32}\right) < \al(N) < \be \cdot \left(1+\frac{1}{32}\right).
$$
In time polynomial in $\log q_M$ we can compute the exact 
expression for $\be$ using rational arithmetic: $\be = p/q$. 
Now we can estimate $N$ and write it as $e^{6 d_k/ 8\be}$ in time polynomial in $\log (q_M)$. 
From there we can continue by adding enough $1$'s to get $I_k$ and 
$\ga_k = [I_k,1,1,\ldots]$. By the construction, it would give us the necessary decrease
 in the value of $r(\ga_k)$. 
\end{proof}

\subsection{Why is $K_\theta$ always computable?}
To provide some intuition why the filled Julia set is computable even when the Julia set is
not, we propose the following toy model.

\begin{figure}[ht]
\begin{center}
\includegraphics[angle=0, scale=0.5]{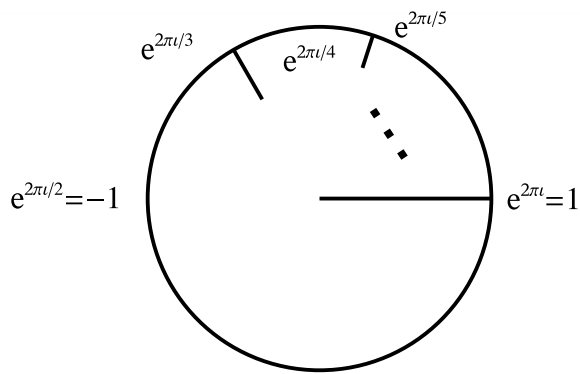}
\end{center}
\caption{Part of the picture of $\Om_t$ for $t=(0.10101\ldots)_2$}
\label{fig:Om}
\end{figure}

Let $A: \NN \ra \{0,1\}$ be any uncomputable 
predicate. Consider the set 
$$\Omega_t=
\left\{
\begin{array}{ll}
S^1\underset{k\in\NN,\;d_k=1}{\bigcup}\{r e^{2\pi i/k}|\; r\in[1-\frac{1}{k},1]\} & \mbox{for } t=(0.d_1d_2d_3\ldots)_2 \in [0,1) \\ \\
S^1\underset{k\in\NN,\;A(k)=1}{\bigcup}\{r e^{2\pi i/k}|\; r\in[1-\frac{1}{k},1]\} & \mbox{for }  t=1 
\end{array}
\right.$$
To avoid ambiguity, we always take the finite expansion for dyadic $t$'s. 
An example of a set $\Om_t$ is depicted on Figure \ref{fig:Om}.
Firstly, note that if $t\in(0,1)$ is not a computable real, then the set
$\Omega_t$ is non-computable by a TM {\it without} an oracle for $t$.
Moreover, even for a TM $M^\phi$ equipped with an oracle input for $t$, the set
$\Om_1$ is clearly non-computable.
 However, when filled, every $\Om_t$ becomes a computable set -- the unit disk. 
This example suggests that filling an uncomputable Julia set
we gain computability at the expense of losing the long and narrow fjords of the Siegel disk.

\subsection{What would a computer really draw?}
One thing is certain -- it is unlikely to draw the true picture of an uncomputable Julia set.
Indeed, by definition of computability, there is no {\it systematic} way of distinguishing
the picture of $J_c$ from all possible pictures on a screen with pixel size $2^{-n}$.
The likelyhood of stumbling upon the correct image when the screen size is, say, 800x600
is very remote. The exact nature of the answer would depend, of course, on the specific
algorithm being used. However, some insight is provided by considering the following
problem, suggested to us by M.~Shub.

Let $\fJ$ be the subset of $\CC \times \CC$ given by
$$
\fJ = \overline{\{ (z,c)~:~z\in J_c\}}. 
$$
Shub has asked us:

\medskip
\noindent
{\it Is the set $\fJ$ computable?}
\medskip

\noindent
The answer again is ``yes'':

\begin{thm}
\label{thm:shub}
Let $d>0$ be any computable real. Then the set
$$\fJ\cap \CC\times \overline{B(0,d)}$$ is a computable subset of $\CC \times \CC$.
\end{thm}

\noindent
Informally, we may think of projection of $\fJ\cap \CC\times(c-\eps,c+\eps)$ 
to the first coordinate as the
picture that a computer could produce when $J_c$ itself is uncomputable.

To understand how the mechanism of non-computability is destroyed in this case,
consider again the set $\Omega_t$ for $t\in(0,1]$ as the toy model.
The set $$W=\overline{\{ (z,t)~:~z\in \Omega_t,\;t\in(0,1]\}}\subset \CC\times\RR$$ is  
computable even though $\Omega_t$ itself is non-computable for $t=1$. This happens
because in the closure of $W$ the ``slice" corresponding to $t=1$ is 
$$
S^1\underset{k\in\NN}{\bigcup}\{r e^{2\pi i/k}|\; r\in[1-\frac{1}{k},1]\} \supset \Om_1.
$$
This set ``masks" the computational hardness of $\Om_1$, and makes $W$ 
computable.

\subsection*{Computability of the set $\fJ$}
We prove \thmref{thm:shub} by showing that $\fJ$ is {\em weakly} computable (Definition \ref{def:weak}). 

We will need the following lemma. 

\begin{lem}
\label{lem:onlyhyp}
For any point $(z,c)$ in the complement of the closure $\overline{\fJ}$, 
$z$ converges to an attracting periodic orbit of $f_c: z \mapsto z^2+c$. 
\end{lem}

The proof of the lemma will require us to recall the nature of discontinuities in the function $\JJ(c)$, particularly the
theory of parabolic implosion. We will give it in the end of the section.

The following lemma allows us to ``cover" all points that belong to
$\fJ$.

\begin{lem}
\label{lem:alg1}
There is an algorithm $A_1 (n)$ that  on input $n$ outputs a sequence of 
dyadic points $p_1, p_2, \ldots \in \CC \times \CC$ such that 
$$
B(\fJ,2^{-(n+3)}) \subset \bigcup_{j=1}^{\infty} B(p_j, 2^{-(n+2)}) \subset B(\fJ, 2^{-(n+1)}). 
$$
\end{lem}

\begin{proof}
The repelling periodic orbits of $f_c$ are dense in $J_c$. 
Hence, the set 
$$
S_{rep} = \{(z,c)~:~z\mbox{ is in a repelling periodic orbit of }f_c\}
$$
is dense in $\fJ$. $S_{rep}$ is a union of a countable number of algebraic 
curves $S^m_{rep}$ given by the constraints 
$$
\left\{
\begin{array}{l}
f_c^m(z)=z \\
|(f_c^m)'(z)|>1
\end{array}
\right.
$$
For each $m$ we can compute a finite number of points $p_1^m, \ldots, p_{r_m}^m$ approximating 
$S^m_{rep}$ such that 
$$
B(S^m_{rep},2^{-(n+3)}) \subset \bigcup_{j=1}^{r_m} B(p_j^m, 2^{-(n+2)}) \subset B(S^m_{rep}, 2^{-(n+1)}). 
$$
We have 
$$
\overline{\fJ} = \overline{S_{rep}} = \overline{\bigcup_{m=1}^\infty S^m_{rep}}. 
$$
Hence the computable sequence $p_1^1, \ldots, p_{r_1}^1, p_1^2, \ldots, p_{r_2}^2, \ldots, p_1^m, \ldots, p_{r_m}^m,\ldots$ satisfies the conditions of the lemma. 
\end{proof}

\begin{cor}
\label{cor:alg1}
There is an oracle machine $M_1^{\phi_1, \phi_2}(n)$, where $\phi_1$ is an oracle for 
$z \in \CC$ and $\phi_2$ is an oracle for $c \in \CC$, such that $M_1^{\phi_1, \phi_2}$ always
halts whenever $d((z,c),\fJ)<2^{-(n+4)}$ and never halts if $d((z,c),\fJ)\ge 2^{-n}$.
\end{cor}

\begin{proof}
Query the oracles for a point $p \in \CC\times \CC$ such that $d(p,(z,c))<2^{-(n+4)}$. 
Then run the following loop:

\medskip
\noindent
$i \leftarrow 0$ \\
 {\bf do}\\
$\mbox{~~~~~~~}$ $i \leftarrow i+1$\\
$\mbox{~~~~~~~}$ generate $p_i$ using $A_1(n)$ from Lemma \ref{lem:alg1} \\
{\bf while} $d(p,p_i)>2^{-(n+2)}$   
\medskip

If $d((z,c),\fJ)<2^{-(n+4)}$, then $d(p, \fJ)<2^{-(n+3)}$, hence by Lemma 
\ref{lem:alg1} there is an $i$ such that $d(p,p_j)\le 2^{-(n+2)}$, and the 
loop terminates. If $d((z,c),\fJ)>2^{-n}$, then $d(p,\fJ)>2^{-n}-2^{-(n-4)} >1.5 \cdot 2^{-(n+1)}$. 
Hence, by Lemma \ref{lem:alg1}, $p \notin B(p_i, 2^{-(n+1)})$ for all $i$, and the loop will 
never terminate. 
\end{proof}

The following lemma allows us to exclude points outside $\overline{\fJ}$ from 
$\fJ$. 

\begin{lem}
\label{alg2}
There is an oracle machine $M_2^{\phi_1, \phi_2}$, where $\phi_1$ is an oracle for 
$z \in \CC$ and $\phi_2$ is an oracle for $c \in \CC$, such that $M_2^{\phi_1, \phi_2}$
halts if and only if $z$ converges to an attracting periodic orbit (or to $\infty$)
under $f_c: z \mapsto z^2+c$. 
\end{lem}

\begin{proof}
$M_2$ is systematically looking for an attracting cycle of $f_c$.
It also iterates $f_c$ on $z$ with increasing precision and for 
increasingly many steps until we are sure that either one of 
the two things holds:
\begin{enumerate}
\item 
the orbit of $z$ converges to $\infty$; or
\item 
we find an attracting orbit of $f_c$ and the orbit of $z$ converges to it. 
\end{enumerate}

If  the search is done systematically, the machine will eventually halt
if one of the possibilities above holds. It obviously won't halt if neither 
holds.
\end{proof}

\begin{proof}[Proof of Theorem \ref{thm:shub}] {\bf The algorithm is:} 
Run the machines $M_1^{\phi_1, \phi_2}(n)$ from Corollary \ref{cor:alg1} and $M_2^{\phi_1, \phi_2}$
from Lemma \ref{alg2} in parallel. Output $1$ if $M_1$ terminates first and $0$ if $M_2$ terminates
first. 

First we observe that $M_1(n)$ only halts on points that are $2^{-n}$-close to $\fJ$, in which 
case $1$ is a valid answer according to Definition \ref{def:weak}. Similarly, $M_2$ only halts
on points that are outside $\fJ$, in which case $0$ is a valid answer. Hence if the 
algorithm terminates, it outputs a valid answer. It remains to see that it does always terminate. 
Consider two cases.

{\bf Case 1: $(z,c) \in \overline{\fJ}$.} In this case $d((z,c),\fJ)=0<2^{-(n+4)}$, and 
the first machine is guaranteed to halt. 

{\bf Case 2: $(z,c) \notin \overline{\fJ}$.} By Lemma \ref{lem:onlyhyp}, $z$ converges to 
an attracting periodic orbit of $f_c$ in this case, and hence the second machine is 
guaranteed to halt. 
\end{proof}

\subsection*{Proof of Lemma \ref{lem:onlyhyp}}

Suppose $z\notin J_c$ and the orbit of $z$ does not belong to an attracting
basin. By the Fatou-Sullivan classification (see e.g. \cite{Mil}), 
there exists $k\in\NN$ such that $w\equiv f_c^k(z)$ belongs to a Siegel
disk or to the  immediate basin of a parabolic orbit.
Our aim is to show that for an arbitrary small $\delta>0$, there exists
a pair $(\tl z,\tl c)\in\CC\times\CC$ with $|z-\tl z|<\delta$,
$|c-\tl c|<\delta$, and for which $\tl z\in J_{\tl c}$.
We will treat the Siegel case first.

\medskip
\noindent
{\it The case when $w$ lies in a Siegel disk}.
Let us denote $\Delta$ the Siegel disk containing $w$, and let $m\in\NN$ be its
period, that is, the mapping $$f_c^m:\Delta\to\Delta$$ is 
conjugated  by a conformal change of coordinates $\phi:\Delta\to \DD$
to an irrational rotation of $\DD$. 

By \propref{siegel discont}, we have the following.
Denote $\zeta=\phi^{-1}(0)\in\Delta$ the center of the Siegel disk.
For each $s>0$ there exists $\tl c\in B(c,s)$ such that
$f_{\tl c}$ has a parabolic periodic point $\tl \zeta$ of period $m$ 
in $B(\zeta,s)$. In particular,
$J_{\tl c}$ is connected, and 
 $B(\zeta,s)\cap J_{\tl c}\neq\emptyset.$

Consider now the $f_c^m$-invariant analytic circle 
$$S_r=\phi^{-1}(\{z=re^{2\pi i\theta},\;\theta\in[0,2\pi)\})$$
which contains $w$. Let $\eps>0$ be such that 
$$B(w,\eps)\subset f_c^k(B(z,\delta))\cap \Delta.$$
Set $B\equiv B(w,\eps/2)$ and 
let $n\in\NN$ be such that the union
$$\bigcup_{0\leq i\leq n}f^{mi}_c(B)\supset S_r.$$
By Proposition \ref{siegel discont} for all $\delta>0$ small enough, there exist
$\tl c\in B(c,\delta)$ for which $J_{\tl c}$ is connected and 
there is a point of $J_{\tl c}$ inside
the domain bounded by $S_r$. Since repelling periodic orbits of $f_c$
are dense in $\partial \Delta$, again for $\delta$ small enough,
there are points of $J_{\tl c}$ on the outside of $S_r$ as well, and
so there exists a point $\xi\in J_{\tl c}\cap S_r$.
By construction, there exists $j\in\NN$ such that 
$f_c^j(B(z,\delta))\ni \xi$. By invariance of Julia set, if
$\tl c$ is close enough to $c$ we have 
$B(z,\delta)\cap J_{\tl c}\neq\emptyset$, and the proof is complete.

\medskip
\noindent
{\it The case when $w$ lies in a parabolic basin.}
Denote $\zeta$ the parabolic periodic point of  
$f_c$ whose immediate basin contains $w$, and let $m\in\NN$ be its
period. We employ the notations of \secref{sec:implosion}.

Recall, that by \thmref{implosion},
for every $s>0$ 
there exists $\tl c\in B(c,s)$ such that
$B(J_{\tl c},s)\supset J_{(c,\tau)}$.

Since $\zeta\in J_c$, and $J_c$ is connected, there exists a point 
$u\in J_c\cap P_R$. Let $\hat w\in C_A$ be the orbit of $w$, and 
let $\hat u\in C_R$ be the orbit of $u$. Choose the translation $\tau:C_A\to C_R$
so that $\tau(\hat w)=\hat u$. Then $J_{(c,\tau)}\ni z$, and the 
claim follows by Theorem \ref{implosion}.

\appendix
\section{Proof of \propref{bdary bded type}}

 Siegel quadratic Julia of bounded type sets may be 
constructed by means of quasiconformal surgery (cf. \cite{Do1}) on a Blaschke product
$$f_\gamma(z)=e^{2\pi i\tau(\gamma)}z^2\frac{z-3}{1-3z}.$$
This map homeomorphically maps the unit circle $\TT$ onto itself with
a single (cubic) critical point at $1$. The angle
$\tau(\gamma)$ can be uniquely selected in such a way that the rotation number
of the restriction $\rho(f_\gamma|_\TT)=\gamma$.

For each $n$, the points $$\{1, f_\gamma(1),f^2_\gamma(1),\ldots,f^{q_{n+1}-1}_\gamma(1)\}$$
form the {\it $n$-th dynamical partition} of the unit circle. The following result is due
to Swiatek and Herman (for the proof see e.g. Theorem 3.1 of \cite{dFdM}):

\begin{thm}[{\bf Universal real {\it a priori} bound}]
\label{real bound}
There exists an explicit constant $B>1$ independent of $\gamma$ and $n$
such that the following holds. 
Let $\gamma\in\RR\setminus\QQ$ and $n\in\NN$. Then
any two adjacent intervals $I$ and $J$
of the $n$-th dynamical partition of $f_\gamma$ are $B$-commensurable:
$$B^{-1}|I|\leq |J|\leq B|I|.$$
\end{thm}

\begin{prop}[\cite{He}]
\label{qs-conjugacy}
For each bounded type $\gamma=[a_0,  \ldots, a_k,\ldots]$ the Blaschke product $f_\gamma$
is $K_1$-quasisymmetrically conjugate to the rotation $R_\gamma:x\mapsto x+\gamma \mod \ZZ$.
The quasisymmetric constant may be taken as $K_1={(2\sup a_i)}^{10B^2}$.
\end{prop}

\noindent
Let us now consider the mapping 
$\Psi$ which identifies the critical orbits of $f_\gamma$ and $P_\gamma$ by
$$\Psi:f^i_\gamma(1)\mapsto P^i_\gamma(c_\gamma).$$

\noindent
We have the following (see, for example, Theorem 3.10 of \cite{YZ}):

\begin{thm}[{\bf Douady, Ghys, Herman, Shishikura}]
\label{surgery}
The mapping $\Psi$ extends to a $K$-quasiconformal homeomorphism of the plane
$\CC$ which maps the unit disk $\DD$ onto the Siegel disk $\Delta_\gamma$.
The constant $K$ may be taken as the quasiconformal dilatation of 
any global quasiconformal extension of the $K_1$-qs conjugacy of \propref{qs-conjugacy}.
In particular, $K\leq 2K_1$.
\end{thm}

Elementary combinatorics implies that each interval of the $n$-th dynamical
partition contains at least two intervals of the $(n+2)$-nd dynamical partition.
This in conjunction with \thmref{real bound} implies that the size of an interval
of the $(n+2)$-nd dynamical partition of $f_\gamma$ is at most $\tau^n$ where
$$\tau=\sqrt\frac{B}{B+1}.$$

Hence, setting
$$\Omega_n=\{P_\gamma^i(c_\gamma),\;i=0,\ldots,q_{n+2}\},$$
by Theorem \ref{surgery}, 
$$ \dist_H ( \Omega_n, \partial \Delta_\gamma) < K \tau^n.  $$

\section{Proof of Lemmas \ref{smlchg}, \ref{notdeclem} and \ref{lem:tailof1s}}
\label{app:B}

We present outlines of proofs for Lemmas \ref{smlchg}, \ref{notdeclem} and \ref{lem:tailof1s}. 
The complete proofs of the intermediate lemmas can be found in \cite{BBY2} and \cite{BY}.
For convenience, we restate the lemmas here:

\smallskip

{\bf Lemma \ref{smlchg}~~}
For any initial segment $I = [a_0, a_1, \ldots, a_n]$, write 
$\omega = [a_0, a_1, \ldots, a_n,$ $ 1, 1, 1, \dots]$. Then 
for any $\ve>0$, there is an $m>0$ and an integer $N$
such that if we write $\be = [a_0, a_1, \ldots, a_n, 1, 1, \ldots, 
1, N, 1, 1, \ldots]$, where the $N$ is located in the $n+m$-th 
position, then
$$
\Phi(\omega) + \ve < \Phi(\be) < \Phi(\omega) + 2 \ve. 
$$

\smallskip

{\bf Lemma \ref{notdeclem}~~}
For  $\omega$ as above, for any $\ve>0$ there is an $m_0>0$, which 
can be computed from $(a_0, a_1, \ldots, a_n)$ and $\ve$,  
such that for any $m \ge m_0$, and for any tail 
$I = [a_{n+m}, a_{n+m+1}, \ldots]$ if we denote 
$$ \be^{I} = [a_1, a_2, \ldots, a_n, 1, 1, \ldots, 1, a_{n+m}, a_{n+m+1}, 
\ldots],$$ then
$$
\Phi(\be^{I}) > \Phi (\om) - \ve. 
$$

\smallskip

{\bf Lemma \ref{lem:tailof1s}~~}
Let $\om = [a_1, a_2, \ldots]$ be a Brjuno number, that is $\Phi(\om)<\infty$. 
Denote $\om_k = [a_1, a_2, \ldots, a_k, 1,1, \ldots]$. 
Then for every $\ve>0$ there is an $m$ such that for all $k\ge m$,
$$
\Phi(\om_k) < \Phi(\om)+\ve. 
$$

\smallskip

Denote
$$
\Phi^{-} (\omega) = \Phi(\omega) - \al_0 (\om) \al_1 (\om) \ldots \al_{n+m-1} (\om) \log 
\frac{1}{\al_{m+n}(\om)}. 
$$
The value of the integer $m>0$ is yet to be determined.
Denote  $$\be^N = (a_0, a_1, \ldots, a_n, 1, 1, \ldots,
1, N, 1, 1, \ldots).$$ 

The following estimates are proven by induction. 

\begin{lem}
\label{4lems}
For any $N$, the following holds:
\begin{enumerate}
\item 
\label{4l:p1}
For $i\le n+m$ we have $$\left| 
\log{\frac{\al_{i}(\be^{N})}{\al_{i}(\be^{N+1})}} \right|
< 2^{i-(n+m)}/N;$$
\item 
\label{4l:p2}
for $i < n+m$,
$$\left| 
\log{\frac{\al_{i}(\be^{N})}{\al_{i}(\be^{1})}} \right|
< 2^{i-(n+m)};$$
\item 
\label{4l:p3}
for $i< n+m$,$$
\left|
\log{\frac{\log \frac{1}{\al_{i}(\be^{N})}}{\log 
\frac{1}{\al_{i}(\be^{N+1})}}} \right|
< 2^{i-(n+m)+1};
$$
\item 
\label{4l:p4}
for $i < n+m-1$, $$
\left|
\log{\frac{\log \frac{1}{\al_{i}(\be^{N})}}{\log 
\frac{1}{\al_{i}(\be^{1})}}} \right|
< 2^{i-(n+m)+1}.
$$
\end{enumerate}
\end{lem}

The estimates yield the following.

\begin{lem}
\label{lem5}
For any $\om$ of the form as in lemma \ref{smlchg} and for any $\ve >0$, 
there is an $m_0 >0$ such that for any $N$ and any $m \ge m_0$,
$$
 | \Phi^{-} (\be^{N}) -
\Phi^{-} ( \be^1) | < \frac{\ve}{4}.
$$
\end{lem}

\begin{proof}
 (Sketch). 
The $\sum$ in the expression for $\Phi (\be^1)$ converges, 
hence there is an $m_1 >1$ such that the tail of the sum 
$\sum_{i\ge n+m_1} \al_1 \al_2 \ldots \al_{i-1} \log \frac{1}{\al_i} < 
\frac{\ve}{16}$. It can be shown that 
\begin{itemize}
\item
for a sufficiently large $m_0>m_1$, if $m>m_0$, then for any $N$ the influence on the 
sum of the ``head" elements is very small:
$$
\left| \sum_{i=1}^{n+m_1-1}  \al_1 (\be^N) \ldots
\al_{i-1}(\be^N) \log
\frac{1}{\al_i(\be^N)} - \sum_{i=1}^{n+m_1-1}  \al_1 (\be^1) \ldots
\al_{i-1}(\be^1) \log
\frac{1}{\al_i(\be^1)} \right| <  
\frac{\ve}{16};
$$
\item 
for the ``tail" terms, for $i\ge n+m_1$ such that $i\neq n+m$,
$$
\frac{\al_1 (\be^N) \ldots \al_{i-1}(\be^N) \log
\frac{1}{\al_i(\be^N)}}{\al_1 (\be^1) \ldots \al_{i-1}(\be^1) \log
\frac{1}{\al_i(\be^1)}} 
 \le e.
$$
\end{itemize}

 After the change each term of the tail could increase 
by a factor of $e$ at most. The value of the ``tail" starts at the 
interval $(0,\frac{\ve}{16}]$, hence it remains 
 in the interval $(0, \frac{e \ve}{16}]$, and the 
change in the tail is bounded  by $\frac{e \ve}{16} < \frac{ 3 \ve}{16}$. 

So the total change in $\Phi^-$ is bounded by
$$
\mbox{change in the ``head"}~~ +~~
\mbox{change in the ``tail"} ~~ < \frac{\ve}{16} + \frac{3 \ve}{16} = 
\frac{\ve}{4}.
$$
\end{proof} 

Lemma \ref{lem5} immediately yields:

\begin{lem}
\label{lem6}
For any $\ve$ and for the same $m_0 (\ve)$ as in lemma \ref{lem5}, 
for any $m \ge m_0$ and $N$, 
$$
 | \Phi^{-} (\be^{N}) -
\Phi^{-} ( \be^{N+1}) | < \frac{\ve}{2}.
$$ 
\end{lem}

Denote $\Phi^{1} (\omega) =  \al_0 (\om) \al_1 (\om) \ldots \al_{n+m-1} (\om) \log 
\frac{1}{\al_{m+n}(\om)} = \Phi(\om) - \Phi^{-} (\om)$. Using the estimates
\ref{4lems} one can prove the following:

\begin{lem}
\label{lem9}
For sufficiently large $m$, for any $N$, 
$$
\Phi^1 (\be^{N+1}) - \Phi^1 (\be^N) < \frac{\ve}{2}.
$$
\end{lem}

Since $\Phi = \Phi^- + \Phi^1$, summing the inequalities  in Lemmas 
\ref{lem6} and \ref{lem9} yields:

\begin{lem}
\label{lem10}
For sufficiently large $m$, for any $N$, 
$$
\Phi (\be^{N+1}) - \Phi (\be^N) < \ve.
$$
\end{lem}

It is immediate from the formula of $\Phi(\beta^N)$ that:

\begin{lem}
\label{lem11}
$$
\lim_{N \rightarrow \infty} \Phi (\be^N) = \infty.$$
\end{lem}

We are now ready to prove Lemma \ref{smlchg}. 

\begin{proof}[Proof of Lemma \ref{smlchg}.] Choose $m$ large enough for lemma 
\ref{lem10} to hold. Increase $N$ by one at a time starting with $N=1$. 
We know that $\Phi(\be^1) = \Phi(\om) < \Phi(\om) + \ve$, and 
by Lemma \ref{lem11}, there exists an $M$ with $\Phi(\be^M) > \Phi(\om) + 
\ve$. Let $N$ be the smallest such $M$.  Then $\Phi(\be^{N-1}) \le
\Phi(\om) +\ve$, and by lemma \ref{lem10}
$$
\Phi(\be^N) < \Phi(\be^{N-1}) + \ve \le \Phi(\om) + 2 \ve.
$$
Hence 
$$
\Phi(\om) +\ve < \Phi(\be^N) <  \Phi(\om) + 2 \ve.
$$
Choosing $\be = \be^N$ completes the proof. 
\end{proof}

\noindent
The second part of the following Lemma follows by the same argument as Lemma \ref{lem5} 
by taking $N\ge 1$ to be an arbitrary real number, not necessairily
an integer. The first part is obvious, since the tail of $\om$ has only $1$'s. 

\begin{lem}
\label{auxbdtl1}
For an $\om = \be^1$ as above, for any $\ve>0$ there is an $m_0 >0$, 
such that for any $m \ge m_0$, and for any tail 
$I = [a_{n+m}, a_{n+m+1}, \ldots]$ if we denote 
$$ \be^{I} = [a_1, a_2, \ldots, a_n, 1, 1, \ldots, 1, a_{n+m}, a_{n+m+1}, 
\ldots],$$ then
$$
\sum_{i \ge n + m} \al_1 (\be^1) \al_2 (\be^1) \ldots  \al_{i-1} (\be^1)
\log \displaystyle\frac{1}{\al_i(\be^1)} < \ve, 
$$
and 
$$
\sum_{i=1}^{n+m-1} \left| \al_1 (\be^I) \ldots
\al_{i-1}(\be^I) \log
\displaystyle\frac{1}{\al_i(\be^I)} - \al_1 (\be^1) \ldots
\al_{i-1}(\be^1) \log
\displaystyle\frac{1}{\al_i(\be^1)} \right| < 
\ve.$$
\end{lem}

We can now prove Lemma \ref{notdeclem}.

\begin{proof}[ Proof of Lemma \ref{notdeclem}.] 
Applying lemma \ref{auxbdtl1} with $\frac{\ve}{2}$ instead of $\ve$, we get
$$
\Phi(\be^I) - \Phi(\om) = \sum \{ \mbox{``head"($\be^I$)}-\mbox{``head"($\om$)}\} + 
\sum \{ \mbox{``tail"($\be^I$)}-\mbox{``tail"($\om$)}\} > 
$$
$$
-\frac{\ve}{2}-\sum \{\mbox{``tail"($\om$)}\} > -\frac{\ve}{2}-\frac{\ve}{2} = -\ve. 
$$
\end{proof}

\begin{proof}[Sketch of the proof of Lemma \ref{lem:tailof1s}.] 
We divide the sum for $\Phi(\om)$,
$$
\Phi(\om)= \underbrace{\sum_{i=1}^s \al_1(\om) \ldots \al_{i-1}(\om) \log \frac{1}{\al_i(\om)}}_{\mbox{``head"}}+
\underbrace{\sum_{i=s+1}^\infty \al_1(\om) \ldots \al_{i-1}(\om) \log \frac{1}{\al_i(\om)}}_{\mbox{``tail"}},
$$
so that $\mbox{``tail"}<\ve/16$. Using the estimates from Lemma \ref{4lems} one can show that 
modifying $\om$ to $\om_k$ for some appropriately chosen $k \gg s$ will satisfy:
\begin{itemize}
\item 
$$\sum_{i=1}^s \al_1(\om_k) \ldots \al_{i-1}(\om_k) \log \frac{1}{\al_i(\om_k)} <
\sum_{i=1}^s \al_1(\om) \ldots \al_{i-1}(\om) \log \frac{1}{\al_i(\om)}+\ve/16,$$
since the relative error in the ``head" terms can be made arbitrarily small;
\item 
for $i>s$,
$$
\al_1(\om_k) \ldots \al_{i-1}(\om_k) \log \frac{1}{\al_i(\om_k)} < 9 \cdot 
\sum_{i=1}^s \al_1(\om) \ldots \al_{i-1}(\om) \log \frac{1}{\al_i(\om)}+ 2^{-k/2}.
$$
Note that for $i>k$ the $2^{-k/2}$ term alone dominates the expression on the left. 
\end{itemize}
Finally for a $k$ as above,
$$
\Phi(\om_k) < \mbox{``head"} (\om_k)+\mbox{``tail"} (\om_k)< \mbox{``head"} (\om)+\ve/16+9 \mbox{``tail"}(\om)+
2^{-s+1} <$$ $$\mbox{``head"} (\om) +\ve/16+ 9 \ve/16+2^{-s+1}<\Phi(\om)+\ve,
$$
for a sufficiently large $s$. 

\end{proof}


\newpage

\begin{thebibliography}{*****}
\bibitem[Ahl]{Ahl} L. Ahlfors. {\it Lectures on Quasiconformal Mappings}, Amer. Math. Soc., 2006.

\bibitem[ABC]{ABC} A. Avila, X. Buff, A. Ch{\'e}ritat. Siegel disks with smooth boundaries. {\it  Acta Math.}  
{\bf 193}(2004),  no. 1, 1--30.

\bibitem[BM]{BM} S. Banach, S. Mazur. Sur les fonctions caluclables. {\it Ann. Polon. Math.} {\bf 16}(1937)

\bibitem[BBY1]{BBY1} I. Binder, M. Braverman, M. Yampolsky. Filled Julia sets with empty interior are
computable. e-print, math.DS/0410580. To appear in {\em Journ. of FoCM}.

\bibitem[BBY2]{BBY2} I. Binder, M. Braverman, M. Yampolsky. 
On computational complexity of Siegel Julia sets. {\it Commun. Math. Phys.},  
{\bf 264}(2006),  no. 2, 317--334.

\bibitem[BB]{BB} E. Bishop, D.S. Bridges. {\it Constructive Analysis}. Springer-Verlag, Berlin (1985)

\bibitem[BCSS]{BCSS} L. Blum, F. Cucker, M. Shub, S. Smale, {\it Complexity and Real Computation}, 
Springer, New York, 1998.

\bibitem[BC1]{BC2} X. Buff, A. Ch{\'e}ritat,  Quadratic Siegel disks with smooth boudaries.
Preprint Univ. Paul Sabatier, Toulouse, III, Num. 242.

\bibitem[BC2]{BC} X. Buff, A. Ch{\'e}ritat, The Brjuno function continuously estimates the size of quadratic Siegel Disks, {\it Annals of Math.} {\bf 164}(2006), No. 1, 265-312.

\bibitem[Bru]{Bru} A.D. Brjuno. Analytic forms of differential equations, {\it Trans. Mosc. Math. Soc} {\bf 25}(1971)

\bibitem[Brv]{thesis}
M. Braverman, ``Computational Complexity of Euclidean Sets: Hyperbolic Julia Sets are 
Poly-Time Computable", M. Sc. Thesis, University of Toronto, 2004, and Proc. CCA 2004,
to appear.

\bibitem[Brv2]{Brv-FOCS} M. Braverman, On  the complexity of real functions, {\it Proc. 46 annual IEEE Symposium
FOCS'05}, 155-164.

\bibitem[Brv3]{Brv2} M. Braverman, 
Parabolic Julia Sets are Polynomial Time Computable. e-print math.DS/0505036. To appear in  {\em Nonlinearity} {\bf 19}, 2006.

\bibitem[BrC]{BrC} M. Braverman, S. Cook.  Computing over the Reals: Foundations for Scientific Computing.  
{\em Notices of the AMS}, {\bf 53}(3),  2006.

\bibitem[BY]{BY} M. Braverman, M. Yampolsky. Non-computable Julia sets. {\it Journ. Amer. Math. Soc.}  {\bf 19}(2006),  no. 3, 551--578 

\bibitem[CK]{CK} A. Chou, Ker-I Ko,  Computational complexity of two-dimensional regions.  {\it SIAM J. Comput.} {\bf  24}(1995),  no. 5, 923--947. 

\bibitem[dFdM]{dFdM} E. de~Faria and W. de~Melo. Rigidity of critical
circle mappings I. {\it J. Eur. Math. Soc. (JEMS)} {\bf 1}(1999), no.    4, 339-392. 



\bibitem[Dou1]{Do1} A. Douady.  Disques de Siegel et anneax de Herman, {\it Sem. Bourbaki, Ast\'erisque}, {\bf 152-153}(1987), 151-172.



\bibitem[Dou2]{Do} A. Douady. Does a Julia set depend continuously on the polynomial? In {\it Complex
dynamical systems: The mathematics behind the Mandelbrot set and Julia sets}. ed. R.L. Devaney,
Proc. of Symposia in Applied Math., Vol 49, Amer. Math. Soc., 1994, pp. 91-138.

\bibitem[DH1]{orsay-notes} A. Douady, J.H. Hubbard. Etude dynamique des polyn\^omes complexes, I-II.
Pub. Math. d'Orsay, 1984.

\bibitem[DH2]{DH2} A. Douady, J.H. Hubbard. On the dynamics of polynomial-like mappings.
{\it Ann. Sci. {\'E}c. Norm. Sup.}, {\bf 18}(1985), 287-343.

\bibitem[Grz]{Grz} A. Grzegorczyk,  Computable functionals, {\it Fund. Math.} {\bf 42}, pp. 168-202, 
1955. 

\bibitem[He]{He} M. Herman. Conjugaison quasi sym{\'e}trique des hom{\'e}omorphismes du cercle {\`a} des rotations, 
Manuscript, 1986. and Conjugaison quasi sym{\'e}trique des diff{\'e}omorphismes du cercle {\`a}
 des rotations et applications aux disques singuliers de Siegel, Manuscript, 1986. Available from \\
{\sl http://www.math.kyoto-u.ac.jp/$\sim$mitsu/Herman/index.html}


\bibitem[Ko1]{Ko} K. Ko, {\it Complexity Theory of Real Functions}, Birkh\"{a}user, Boston, 1991.

\bibitem[Ko2]{Kosurv}
K. Ko, {Polynomial-time computability in analysis}, in 
''Handbook of Recursive
        Mathematics", Volume {\bf 2} (1998), Recursive Algebra, Analysis and 
Combinatorics, Yu. L. Ershov et al.
        (Editors), pp 1271-1317. 


\bibitem[Lac]{Lac}
D. Lacombe, Extension de la notion de fonction r{\'e}cursive aux fonctions d'une
ou plusiers variables, {\it C. R. Acad. Sci. Paris}, {\bf 240}, pp. 2473-2480, {\bf 241}, 13-14 and 151-153, 1955. 

 

\bibitem[MMY]{MMY} S. Marmi, P. Moussa, J.-C. Yoccoz, The Brjuno functions 
and their regularity properties, {\it Commun. Math. Phys.}
{\bf 186}(1997), 265-293.

\bibitem[Mat]{Mat} Y. Matiyasevich, {\it Hilbert's Tenth Problem}, The MIT Press, Cambridge, London, 
1993. 

\bibitem[Maz]{Maz} S. Mazur, {\it Computable analysis}, Rosprawy Matematyczne, Warsaw, vol. 33 (1963).



\bibitem[McM1]{McM1} C. McMullen. {\it Complex dynamics and renormalization}. {\it Annals
of Math. Studies}, v.135, Princeton Univ. Press, 1994.

\bibitem[Mil]{Mil} J. Milnor. {\it Dynamics in one complex variable. Introductory lectures.} 
Princeton University Press, 2006.

\bibitem[Pap]{Papad} C. H. Papadimitriou, {\it Computational Complexity}, Addison Wesley, 1994.

\bibitem[Pet]{Pet} C. Petersen,  Local connectivity of some Julia sets 
containing a circle with an irrational rotation, {\it Acta Math.}, {\bf 177} (1996) 163-224.

\bibitem[PZ]{PZ} C. Petersen, S. Zakeri.  On the Julia set of a typical quadratic polynomial with a Siegel disk.  
{\it Ann. of Math.} (2) {\bf 159}(2004), no. 1, 1--52.

\bibitem[Pom]{Pom} C. Pommerenke, {\it Boundary behaviour of conformal maps}, Springer-Verlag, 1992.

\bibitem[RW]{WeiPaper}
R. Rettinger, K. Weihrauch,  The Computational Complexity of Some Julia 
Sets, in {\it  ACM STOC'03}, San Diego, CA, 2003.


\bibitem[Ret]{Ret} R. Rettinger.  A fast algorithm for Julia sets of hyperbolic rational functions, {\it Proc. of CCA 2004},  145--157,
Electron. Notes Theor. Comput. Sci., 120,
Elsevier, Amsterdam, 2005.

\bibitem[RZ]{RZ} S. Rohde, M. Zinsmeister.  Variation of the conformal radius.  {\it J. Anal. Math.}  {\bf 92}(2004), 105--115.

\bibitem[Sie]{siegel} C. Siegel, Iteration of analytic functions.  {\it Ann. of Math.} (2)  43,  (1942). 607--612

\bibitem[Sip]{Sip} M. Sipser, {\it Introduction to the Theory of Computation}, PWS Publishing 
Company, 1997.

\bibitem[Shi]{Shi} M. Shishikura, Bifurcation of parabolic fixed points.  {\it The Mandelbrot set, theme and variations,}  325--363, London Math. Soc. Lecture Note Ser., 274, Cambridge Univ. Press, Cambridge, 2000.

\bibitem[S{\o}r]{Sor} D. S{\o}rensen. Describing quadratic Cremer point polynomials by parabolic perturbations. {\it  Ergodic Theory Dynam. Systems}  {\bf 18}(1998),  no. 3, 739--758. 

\bibitem[Sul]{Sul} D. Sullivan, Conformal dynamical systems, in {\it Geometric Dynamics,} 
Ed. Palis, Lecture Notes Math., {\bf 1007}(1983), Springer-Verlag, 725-752.

\bibitem[Tur]{Tur} A. M. Turing, On Computable Numbers, With an Application to the 
Entscheidungsproblem. In {\it Proceedings, London Mathematical Society}, 1936, pp. 230-265.

\bibitem[Wei]{Wei} K. Weihrauch, {\it Computable Analysis}, Springer, Berlin, 2000.

\bibitem[Wey]{Wey} H. Weyl. {\it Randbemerkungen zu Hauptproblemen der Mathematik, II, Fundamentalsatz der Algebra and
Grundlagen der Mathematik}, Math. Z., {\bf 20}(1924), pp. 131-151.

\bibitem[Yam]{Yam} M. Yampolsky.
Complex bounds for renormalization of critical circle maps,
{\it Erg. Th. \&  Dyn. Systems}. {\bf 19}(1999), 227-257.



\bibitem[Yoc]{Yoc} J.-C. Yoccoz, {\it Petits diviseurs en dimension 1}, {S.M.F., Ast{\'e}risque}, {\bf 231}(1995).

\bibitem[YZ]{YZ} Yampolsky, Zakeri,  Mating Siegel quadratic polynomials.
{\it J. Amer. Math. Soc.} {\bf 14}(2001), no. 1, 25--78 

\end{thebibliography}
\end{document}